\documentclass[10pt]{amsart}
\usepackage{amsthm}
\usepackage{amsmath}
\usepackage{amssymb}
\usepackage{pstricks, pst-node}

\usepackage{graphicx}

\theoremstyle{plain}
 \newtheorem{theorem}{Theorem}[section]
 
  \newtheorem{proposition}[theorem]{Proposition}
 
 \newtheorem{conjecture}[theorem]{Conjecture}
 \newtheorem{openproblem}[theorem]{Open Problem}
 \theoremstyle{definition}
  \newtheorem{definition}[theorem]{Definition}
  \newtheorem{example}[theorem]{Example}
  
 \theoremstyle{remark}
  \newtheorem{remark}[theorem]{Remark}

\newcommand{\R}{\mathbb{R}}
\newcommand{\N}{\mathbb{N}}
\newcommand{\PP}{\mathbb{P}}
\newcommand{\ZZ}{\mathbb{Z}}
\newcommand{\CC}{\mathbb{C}}

\newcommand{\VV}{\mathcal{V}}

\newcommand{\xx}{\mathbf{x}}

\newcommand{\Sym}{\mathfrak{S}}
\newcommand{\lieg}{\mathfrak{g}}
\newcommand{\nn}{\mathfrak{n}}

\newcommand{\diag}{\mathrm{diag}}
\newcommand{\character}{\mathrm{char}}

\newcommand{\vicin}{\prec}

\newcommand{\vertex}{\mathrm{vert}}
\newcommand{\subfacet}{\mathrm{subf}}

\begin{document}

%\date{August 2006}

\title[Shifted families, degree sequences, and plethysm]
        {Shifted set families, degree sequences, and plethysm}

\author{C. Klivans}
\address{Depts. of Mathematics and Computer Science \\Univ. of Chicago}

\author{V. Reiner}
\address{School of Mathematics\\Univ. of Minnesota}

\thanks{First author partially supported by NSF VIGRE grant
DMS-0502215. Second author partially supported by NSF grant
DMS-0601010.}


\begin{abstract}
  We study, in three parts, degree sequences of $k$-families 
(or $k$-uniform hypergraphs) and shifted $k$-families.  
\begin{enumerate}
\item[$\bullet$] 
The first part collects for the first time in one place, various implications 
such as 
$$ 
\textrm{Threshold} \Rightarrow  \textrm{Uniquely Realizable}
\Rightarrow  \textrm{Degree-Maximal} \Rightarrow  \textrm{Shifted} 
$$
which are equivalent concepts for $2$-families (= simple graphs), but 
strict implications for $k$-families with $k \geq 3$.  The
implication that uniquely realizable implies degree-maximal seems to be new.
\item[$\bullet$]
The second part recalls Merris and Roby's reformulation of the characterization
due to Ruch and Gutman for graphical degree sequences and shifted $2$-families.
It then introduces two generalizations which are characterizations
of shifted $k$-families.

\item[$\bullet$] 
The third part recalls the connection between degree sequences of $k$-families
of size $m$ and the plethysm of elementary symmetric functions $e_m[e_k]$.
It then uses highest weight theory to explain how
shifted $k$-families provide the ``top part'' of these plethysm expansions,
along with offering a conjecture about a further relation.
\end{enumerate}
\end{abstract}

\maketitle

\tableofcontents

\section{Introduction}

Vertex-degree sequences achievable by simple graphs are well-understood and characterized--
e.g. \cite{Sierksma} offers seven equivalent characterizations.
By contrast, vertex-degree sequences achievable by simple {\it hypergraphs} are poorly understood, 
even for $k$-uniform hypergraphs ($k$-families), even for $k=3$. 

The current paper has three goals/parts.  The first part is about
various equivalent concepts for graphs such as {\it positive threshold, 
threshold, 
uniquely realizable, 
degree-maximal}, and 
{\it shifted}
which arise in the literature as 
the extreme cases in characterizations of degree sequences.
Here our goal (Theorem~\ref{omnibus-theorem}) is to explain how these turn into  
a strict hierarchy of concepts for $k$-families when $k>2$.
Most of the implications in the hierarchy have occurred in scattered
places in the literature, although one of them (uniquely realizable implies
degree-maximal) appears to be new.  After defining the relevant concepts in
Section~\ref{definition-section}, Theorem~\ref{omnibus-theorem} is proven in
Section~\ref{relations-section}.

The second part (Section~\ref{characterization-section}) addresses
 characterizing degree sequences for $k$-families more explicitly and
makes a promising start on this
problem.  Proposition~\ref{unsatisfactory-proposition} offers a
reduction to shifted families stating that an integer sequence is
a degree sequence if and only if it is majorized by a shifted degree
sequence. Such shifted sequences are unfortunately also poorly
understood. This section then re-examines Merris and Roby's
reformulation of Ruch and Gutman's characterization of graphical
($k=2$) degree sequences, as well as their characterization of the
extreme case of shifted graphs.  Given an integer partition, Merris
and Roby's conditions are stated in terms of the associated Ferrers diagram.
The goal in this part is to prove the more general
Proposition~\ref{k-family-geometry}, giving a $k$-dimensional
extension for shifted $k$-families via associated stacks of cubes.

The third part (Section~\ref{plethysm-section}) 
recalls a related and well-known connection between graph degree sequences
and the $k=2$ case of the problem of 
expanding {\it plethysms} $e_m[e_k]$ of elementary symmetric functions
in terms of Schur functions $s_\lambda$.  This problem was solved 
by a famous identity due to Littlewood:
$$
\sum_{\substack{\text{ all simple} \\\text{ graphs }K}} \xx^{d(K)} 
    \left( = \prod_{i<j} (1 + x_i x_j) = \sum_{m \geq 0} e_m[e_2] \right)
    = \sum_{\substack{\text{ shifted} \\\text{ graphs }K}}  s_{d(K)}.
$$
The goal in the third part is to prove that the generalizations
for $k>2$ of the left and right sides in this identity,
$$
\sum_{\substack{\text{ all }k-\text{uniform} \\\text{hypergraphs }K}} \xx^{d(K)} 
  \left( = \prod_{\substack{k-\text{subsets} \\\{i_1, i_2, \cdots, i_k\}}} 
        \left( 1 + x_{i_1} x_{i_2} \cdots x_{i_k} \right) 
    = \sum_{m \geq 0} e_m[e_k] \right)
$$
and
$$
\sum_{\substack{\text{ shifted }k-\text{uniform} \\\text{ hypergraphs }K}}  s_{d(K)},
$$
while {\it not} being equal, do have many properties in common.  In particular,
they 
\begin{enumerate}
\item[$\bullet$] have the same monomial support (Proposition~\ref{small-k-coincidences}), 
\item[$\bullet$] both enjoy two extra symmetries 
(Propositions~\ref{symmetry-one} and \ref{symmetry-two}), 
\item[$\bullet$] have the Schur expansion for the former 
coefficientwise larger than for the latter (Theorem~\ref{shifteds-as-lower-bound}).
\end{enumerate}

\section{Definitions and Preliminaries}
\label{definition-section}

\subsection{The basic definitions}

After defining $k$-families and degree sequences, we
recall some of the basic definitions.

\begin{definition} ($k$-families)
Let $\PP:=\{1,2,\ldots\}$ and $[n]:=\{1,2,\ldots,n\}$.
A $k$-family $K$ on $[n]$ is a collection
$K=\{S_1,\ldots,S_m\}$ of distinct $k$-subsets  $S_i \subset [n]$.
In other words $S_i \in \binom{[n]}{k}$.
These are sometimes called {\it (simple) $k$-uniform hypergraphs},
and the $S_i$ are called the {\it hyperedges}.
Say that $K$ has {\it size} $m$ if $|K|=m$.

Two $k$-families $K, K'$ are {\it isomorphic} if there exists a permutation
$\sigma$ of $[n]$ which relabels one as the other: $\sigma(K)=K'$.
\end{definition}

\begin{definition} (Degree sequence) \rm \
For a simple graph $G = (V,E)$ with $|V| = n$,  the vertex-degree
sequence of $G$ is the sequence $d(G) = (d_1, d_2, \ldots, d_n)$ where
$d_i = |\{j : \{i,j\} \in E \}|$.  
More generally, the {\it (vertex-) degree sequence} for a $k$-family $K$ on $[n]$ is
$$d(K) = (d_1(K), d_2(K), \ldots, d_n(K))$$
where $d_i(K) = |\{S \in K : i \in S \}|.$

For any integer sequence $d=(d_1,\ldots,d_n)$, let $|d|:=\sum_{i=1}^n d_i$
denote its sum or {\it weight}.

\end{definition}

With these definitions in hand, we define the main conditions on
$k$-families to be studied here.

\begin{definition}  (Threshold families) \rm \
Given a $k$-subset $S$ of $[n]$, its {\it characteristic vector}
$\chi_S \in \{0,1\}^n$ is the sum of standard basis vectors 
$\sum_{i \in S} e_i$.  In other words, $\chi_S$ is the vector of length $n$ 
with ones in the coordinates indexed by $S$ and zeroes in all other coordinates.  
Note that $d(K) = \sum_{S \in K} \chi_S.$

  A $k$-family $K$ of $[n]$ is {\it threshold} if there exists a linear functional
  $w \in (\R^n)^*$ such that $S \in K$ if and only if $w(\chi_S) > 0$.  

  A variation on this was introduced by Golumbic \cite[Property $T_1$,
  page 233]{Golumbic} and studied by Reiterman, R\"{o}dl,
  \v{S}i\v{n}ajov\'{a}, and T\.{u}ma \cite{Reiterman}.  Say that $K$
  is {\it positive threshold} if there is a linear functional
  $w(x)=\sum_{i=1}^n c_i x_i$ having {\it positive} coefficients $c_i$
  and a {\it positive} real threshold value $t$ so that $S \in K$ if
  and only if $w(\chi_S) > t$.

\end{definition}

\begin{example} \rm \
Consider a $k$-family of $[n]$ that consists of all possible
$k$-sets.  Such ``complete'' families are threshold: simply take
any strictly positive linear functional.  The empty family is
similarly threshold, as can be seen by taking any strictly negative linear functional.

The $3$-family $\{123, 124, 125 \}$ is threshold with $w =
(1, 1, -1, -1, -1)$.  This example may be extended to general $k$ by
taking a family of $k$-sets which have a common $(k-1)$-set in their
intersection.  For this family take the linear functional that weights the vertices in the
common $(k-1)$-set with $1$ and all other vertices with $-(k-2)$.

\end{example}

\begin{definition} (Uniquely realizable families) \rm \ 
A $k$-family $K$ is  uniquely realizable if there does not exist a $k$-family $K' \neq K$
with $d(K) = d(K')$. 
\end{definition}

\begin{example} \rm \
It is possible to have two non-isomorphic families with the same degree
sequence.  Let $K$ be a
disjoint union of two cycles of length $3$ and $K'$ be a cycle of
length $6$.  Both families have degree sequence $(2,2,2,2,2,2)$ and hence
 are not uniquely realizable.

It is not necessary, however, to consider non-isomorphic families.  The
$2$-family $K = \{12,23,34\}$, a path of length $3$, with degree
sequence $(1,2,2,1)$ is not uniquely realizable. The $2$-family $K' =
\{13,23,24\}$, also a path of length three, has the same degree sequence.

The family $K = \{12,23,13\}$, a single cycle of length $3$, which has degree sequence $(2,2,2)$ is uniquely realizable.
\end{example}

Note that two $k$-families $K$ and $K'$ of the same size $m=|K| =
|K'|$ will have the same sum for their degree sequences:
$|d(K)|=|d(K')|=km$.  This leads naturally to considering the
majorization order for comparing
degree sequences.  Majorization is also known as the dominance order.

\begin{definition} (Degree-maximal families)  \rm \
Given two sequences of real numbers 
$$
\begin{aligned}
a &= (a_1 ,\ldots, a_n),\\
b &= (b_1,\ldots,b_m)
\end{aligned}
$$
with the same sum $|a|=|b|$, one says that $a$ {\em majorizes} $b$ ($a \unrhd b$) 
if the following system of inequalities hold:
$$
\begin{aligned}
a_1 &\geq b_1 \\
a_1 + a_2 &\geq b_1 + b_2\\
&\vdots\\
a_1 + a_2 + \ldots + a_{n-1} &\geq b_1 + b_2 + \ldots + b_{n-1}.
\end{aligned}
$$
\noindent
Write $a \rhd b$ when $a \unrhd b$ but $a \neq b$.

If one weakens the equality $|a|=|b|$ of the total sums to an inequality ($|a| \geq |b|$) 
then one says that $a$ {\em weakly majorizes} $b$ (written $a \succ b$).

A $k$-family $K$ is {\it degree-maximal} if 
there does not exist $K' \neq K$
such that  $d(K') \rhd d(K)$, i.e.  $d(K)$ is maximal with
respect to majorization.
\end{definition}

\begin{example} \rm \
Let $K$ and $K'$ be the $3$-families $\{124, 125, 135\}$ and $\{123,
124, 125\}$, with degree sequences $(3,2,2,1,1)$ and $(3,3,1,1,1)$.
Clearly $d(K') \rhd d(K)$, hence $K$ is not degree-maximal.  It is not
hard to check that $K'$ is degree-maximal.
\end{example}

An important property of the majorization order is that the
weakly decreasing rearrangement of any sequence always majorizes
the original sequence.  A consequence is
that a degree-maximal family $K$ must always have
its degree sequence $d(K)$ weakly decreasing, otherwise the isomorphic
family $K'$ obtained by relabeling the vertices in weakly decreasing order
of degree would have $d(K') \rhd d(K)$.

\begin{definition} (Shifted families) \rm \
The {\it componentwise partial order} (or {\it Gale order}) on the set $\binom{\PP}{k}$
of all $k$-subsets of positive integers is defined as follows: say $x \leq y$ if
$$
\begin{aligned}
x &= \{x_1 < x_2 < \cdots < x_k\},\text{ and }\\
y &= \{y_1 < y_2 < \cdots < y_k\}
\end{aligned}
$$
satisfy $x_i \leq y_i$ for all $i$.

A $k$-family is {\it shifted} if its $k$-sets, when written as increasing strings,
form an order ideal in the componentwise partial order.

When exhibiting a shifted family $K$, if  $\{S_1, S_2, \ldots, S_p \}$ is 
the unique antichain of componentwise maximal $k$-sets in  $K$, 
we will say that $K$ is the shifted family {\it generated by} $\{S_1, S_2, \ldots, S_p \}$,
and write $K=\langle S_1, S_2, \ldots, S_p \rangle$.
\end{definition}

\begin{example} \rm \
\label{first-shifted-example}
The shifted family $K = \langle 235, 146 \rangle$ consists of triples 
$$
\{123, 124, 125, 126, 134, 135, 136, 145, 146, 234, 235\}
$$ 
and has degree sequence $d(K) = (9,6,6,5,4,3)$. 

The family $K = \{124, 125, 134, 135, 234, 235\}$, consisting of the
triples indexing maximal faces in the boundary of a triangular bipyramid, 
is not shifted. The triple $123$ is ``missing'' from the
family.  Furthermore, it is not possible to relabel this family and achieve a
shifted family.
\end{example}

It is an easy exercise to check that a shifted family $K$ will always have
its degree sequence $d(K)$ weakly decreasing.

\subsection{Cancellation conditions}

Here we introduce two cancellation conditions on $k$-families,
which arise in the theory of simple games and {\em weighted} games~\cite{Taylor}.  
Both will turn out to be equivalent
to some of the previous definitions; see Theorem~\ref{omnibus-theorem} below.

\begin{definition} (Cancellation conditions) \rm \

Consider two $t$-tuples of $k$-sets $(A_1,A_2, \ldots, A_t), (B_1,B_2,
\ldots, B_t)$, allowing repetitions in either $t$-tuple, such that
$\sum_{i=1}^t \chi_{A_i} = \sum_{i=1}^t \chi_{B_i}$.
A $k$-family $K$ of $[n]$ satisfies the {\it cancellation condition
$CC_t$} if for any two such $t$-tuples, whenever each $A_j$ is in $K$ then at
least one $B_j$ must also be in $K$.

A $k$-family $K$ satisfies the {\it cancellation condition $DCC_t$}
if for any two collections of $t$ {\it distinct} $k$-sets $\{A_1, \ldots, A_t\},
\{B_1,\ldots,B_t\}$ with $\sum_{i=1}^t \chi_{A_i} = \sum_{i=1}^t \chi_{B_i}$,
whenever each $A_j$ is in $K$ then at least one $B_j$ must also be in $K$.
In the simple games literature this is known as {\it Chow trade-robustness}.
\end{definition}

Note that every $k$-family satisfies $DCC_1(=CC_1)$.  We recall here
the ``simplest'' failures for $DCC_2$, which appear
under the name of {\it forbidden configurations} 
in the study of Reiterman, et al.\cite[Definition 2.3]{Reiterman}.

\begin{definition}
Say that a $k$-family $K$ satisfies the $RRST$-condition if there does not
exist two $(k-1)$-sets $A,B$ and a pair $i,j$ satisfying
$$
\begin{aligned}
&i,j \not\in A\\
&i,j \not\in B\\
&A \sqcup \{j\}, B \sqcup \{i\} \in K, \text{ but }\\
&A \sqcup \{i\}, B \sqcup \{j\} \not\in K.
\end{aligned}
$$

Note that such a tuple $(A,B,i,j)$ would lead to a violation of $DCC_2$ since
$$
\chi_{A \sqcup \{j\}} + \chi_{B \sqcup \{i\}} =
\chi_{A \sqcup \{i\}} + \chi_{B \sqcup \{j\}}.
$$
\end{definition}

\begin{example} \rm \
It is not hard to check that the $3$-family $K = \{123, 134, 145\}$
satisfies $CC_3$ and $DCC_3$.  $K$ does not however satisfy $DCC_2$
as seen by taking the collections $\{123, 145\}$ and $\{135, 124\}$.
\end{example}

\subsection{Vicinal preorder}
In \cite{Klivans} it was shown how shiftedness relates to a certain preorder
on $[n]$ naturally associated to any $k$-family on $[n]$;  see Theorem~\ref{omnibus-theorem}
below.

\begin{definition} (Vicinal preorder) \rm \
\label{vicinal-definition}
Given a $k$-family $K$ on $[n]$ and $i \in [n]$, define the {\it open and closed neighborhoods} 
of $i$ in $K$ to be the following two subcollections $N_K(i), N_K[i]$ of $\binom{[n]}{k-1}$:
$$
\begin{aligned}
N_K(i) & := \left\{ A \in \binom{[n]}{k-1}: A \sqcup \{i\} \in K \right\} \\
N_K[i] & := N_K(i) \sqcup 
   \left\{ A \in \binom{[n]}{k-1}: i \in A \text{ and } A \sqcup \{j\} \in K \text{ for some }j \right\}.
\end{aligned}
$$

Define a binary relation $\vicin_K$ on $[n] \times [n]$ by $i \vicin_K j$ if $N_K[i] \supseteq N_K(j)$.
\end{definition}

\begin{proposition}(\cite[\S 4]{Klivans})
The relation $\vicin_K$ defines a preorder on $[n]$, that is, it is reflexive and transitive.
\end{proposition}
\begin{proof}
Since $N_K[i] \supseteq N_K(i)$, the relation $\vicin_K$ is clearly reflexive.
To show transitivity, assume $N_K[i] \supseteq N_K(j)$ and $N_K[j] \supseteq N_K(k)$, then
we must show $N_K[i] \supseteq N_K(k)$.  Equivalently, we must show that
$$
N_K(k) \cap \left( N_K[j] \setminus N_K(j) \right) \subset N_K[i].
$$
The typical set in $N_K(k) \cap \left( N_K[j] \setminus N_K(j) \right)$ is of the form
$A \sqcup \{j\}$ where $A$ is a $(k-2)$-set for which $A \sqcup \{j,k\} \in K$.
We must show such a set $A \sqcup \{j\}$ lies in $N_K[i]$.

\noindent
{\bf Case 1.} $i \in A$.  

  Then the fact that $\left( A \sqcup \{j\} \right) \sqcup \{k\} \in K$ tells us $A \sqcup \{j\} \in N_K[i]$,
and we're done.

\noindent
{\bf Case 2.} $i \not\in A$.

  Then $A \sqcup \{k\} \in N_K(j) \subseteq N_K[i]$. But $i \not \in A$, so this
forces $A \sqcup \{k\} \in N_K(i)$.  This then implies 
$A \sqcup \{i\} \in N_K(k) \subseteq N_K[j]$.  Since $j \not\in A$, this 
forces $A \sqcup \{i\} \in N_K(j)$.  Hence $A \sqcup \{j\} \in N_K(i) \subseteq N_K[i]$,
as desired.
\end{proof}

\begin{example} \rm \
The shifted family from Example~\ref{first-shifted-example}
$$
K = \langle 235, 146 \rangle = 
\{123, 124, 125, 126, 134, 135, 136, 145, 146, 234, 235\}
$$ 
has its vicinal preorder on $\{1,2,3,4,5\}$ given by
$$
6 \vicin_K 5 \vicin_K 4 \vicin_K 3 \sim_K 2 \vicin_K 1
$$
where we write $i \sim_K j$ if $i \vicin_K j$ and $j \vicin_K i$.  
Note that in this case, the vicinal preorder is a {\it linear
preorder}, that is, every pair of elements $i, j$ are related,
either by $i \vicin_K j$ or by $j \vicin_K i$ or by both.
%Consider the shifted family $K = \langle 125, 234 \rangle$.  
%The vicinal preorder for this family is linear,
%$$ 5 \vicin_K 4 \sim_K 3 \vicin_K 2 \sim_K 1,$$

\end{example} \rm \

\subsection{The zonotope of degree sequences}

Here we recall a zonotope often associated with degree sequences.
For basic facts about zonotopes, see \cite{McMullen}.

\begin{definition} (Polytope of degree sequences)
The {\em polytope of degree sequences} $D_n(k)$ is the convex hull in
$\R^n$ of all degree sequences of $k$-families of $[n]$.
Equivalently, $D_n(k)$ is the zonotope given by the Minkowski sum of
line segments $\{[0,\chi_S] \, | \, S \in \binom{[n]}{k}\}$, where recall
that $\chi_S$ was the sum of the standard basis vectors,  $ \chi_S = \sum_{i \in S} e_i$.

The case $k=2$ was
first considered in~\cite{Koren} and further developed in~\cite{Peled}
and~\cite{Stanley-zonotope}.  The case $k > 2$ was studied more recently in~\cite{Murthy}.
\end{definition}

\subsection{Swinging and shifting}

Certain ``shifting'' operations produce a shifted family from an arbitrary family.  There are two main
variants of shifting: {\it combinatorial shifting}  introduced by
Erd\"{o}s, Ko, and Rado~\cite{EKR} and Kleitman~\cite{Kleitman} and {\it algebraic
shifting} introduced by Kalai~\cite{Kalai}.  Here we consider the
related operation of swinging.

\begin{definition} (Swinging) \rm \ \\
Given a $k$-family $K$ on $[n]$, suppose that there
is a pair of indices $i < j$ and a $(k-1)$-set $A$ containing neither
of $i,j$, such that $A \sqcup \{j\} \in K$ and $A \sqcup \{i\} \notin K$.
Then form the new $k$-family
$$
K' = (K \, \, \backslash \, \, (A \sqcup j)) \cup (A \sqcup i).
$$
In this situation, say that $K'$ was formed by a {\it swing} from $K$.
\end{definition}

The difference between this operation and combinatorial
shifting is the {\it fixed} $(k-1)$-set $A$;  
combinatorial shifting instead chooses a pair of indices $ i < j $ and applies
the swinging construction successively to {\it all} applicable $(k-1)$-sets $A$.
Hence combinatorial shifting is more restrictive: it is
not hard to exhibit examples where a $k$-family $K$ can be associated
with a shifted family $K'$ via a sequence of swings, but not via
combinatorial shifting.  Neither swinging nor combinatorial shifting
is equivalent to algebraic shifting.
Recently, Hibi and Murai \cite{Hibi} have exhibited an example of a
family where the algebraic shift cannot be achieved by combinatorial
shifting. We do not know if all outcomes of algebraic shifting may be
obtained via swinging.

\begin{example} \rm \
Let $K$ be the non-shifted $3$-family $\{123, 124, 145, 156\}$.  First
consider combinatorially shifting $K$ with respect to the pair
$(2,5)$.  The resulting family is
$\{123, 124, 125, 126\}$ and is easily seen to be shifted.

The following swinging operations on $K$ result in a different shifted
family.  First swing $145$ with respect to $(2,4)$ which replaces
$145$ with $125$.  Next swing $156$ with respect to $(3,5)$ which
replaces $156$ with $136$. Finally, swing the new face $136$ with
respect to $(4,6)$.  The result is the shifted family $\{123, 124,
134, 125\}$.

\end{example}

We note here a few easy properties of swinging.

\begin{proposition}
\label{prop:swing}
Assume that the $k$-family $K$ on $[n]$ has been
labelled so that $d(K)$ is weakly decreasing.

\begin{enumerate}
\item[(i)] One can swing from $K$ if and only if $K$ is not shifted.

\item[(ii)] If one can swing from $K$ to $K'$ then $d(K') \rhd d(K)$.  

\item[(iii)]  (\cite[Proposition 9.1]{Duval})
If $d'$ is a weakly decreasing sequence of positive
integers with $d(K) \unrhd d'$, then there exists a family $K'$
with $d(K')=d'$ such that $K$ can be obtained  from $K'$ by a (possibly empty) sequence of swings.
\end{enumerate}
\end{proposition}
\begin{proof}
Assertions (i) and (ii) are straightforward.  We repeat here the proof
of (iii) from \cite[Proposition 9.1]{Duval} for completeness.  Without
loss of generality $d'$ {\it covers} $d(K)$ in the majorization (or
{\it dominance} order) on partitions, which is well-known to imply
\cite{Muirhead} that there exist indices $i < j$ for which
$$
\begin{aligned}
d_i(K)&=d'_i+1,\\
d_{j}(K)&=d'_{j}-1, \\
d_l(K)&=d'_l \text{ for }l \neq i,j.
\end{aligned}
$$
This implies $d_i(K) > d_{j}(K)$, so there must exist at least one $(k-1)$-subset
$A$ for which $A \sqcup \{i\} \in K$ but $A \sqcup \{j\} \not\in K$.
Then perform the reverse swing to produce 
$K':= K \setminus \{A \sqcup \{i\} \} \cup \{ A \sqcup \{j\} \}$ 
achieving $d(K')=d'$.
\end{proof}

\section{Some relations between the concepts}
\label{relations-section}

\begin{theorem} 
\label{omnibus-theorem}
For a $k$-family $K$, the following equivalences and implications hold:
\begin{eqnarray} 
\label{positive-threshold}
          &               & K \text{ is positive threshold}  \\
\label{threshold}
          &\Rightarrow    & K \text{ is threshold}  \\
\label{vertex}
          &\Leftrightarrow& d(K) \text{ is a vertex of }D_n(k)\\
\label{CC}
          &\Leftrightarrow&  K \text{ satisfies }CC_t \text{ for all }t  \\
\label{uniquely-realizable}
          &\Rightarrow    & d(K) \text{ is uniquely realizable} \\
\label{DCC}
          &\Leftrightarrow& K  \text{ satisfies }DCC_t \text{ for all }t   \\
\label{degree-maximal}
          &\Rightarrow    &  K \text{ is isomorphic to 
                                      a degree-maximal family}   \\
\label{vicinal}
          &\Rightarrow    &  K \text{ has its vicinal preorder }\vicin_K\text{ a total preorder}  \\
\label{RRST} 
          &\Leftrightarrow&  K \text{ satisfies }RRST\\
\label{shifted}
          &\Leftrightarrow& K  \text{  is isomorphic to a shifted family}   
\end{eqnarray}
For $k \geq 3$, the four implications shown are strict, while for
$k=2$ these concepts are all equivalent.
\end{theorem}

\begin{remark} \rm \ 
Before proving the theorem, we give references for most of its
assertions.  Only the implication \eqref{uniquely-realizable} (or
equivalently, \eqref{DCC}) implies \eqref{degree-maximal} is new, as
far as we are aware.  Our intent is to collect the above properties
and implications, arising in various contexts in the literature,
together for the first time.

For $k=2$ these concepts describe the class of graphs usually known as
{\it threshold graphs}.  The equivalence of the threshold and shifted
properties for graphs seems to have been first noted
in~\cite{Chvatal}. Properties $(2), (3), (5), (7), \, \textrm{and} \,
(8)$ for graphs may be found in~\cite{Mahadev}.  Properties $(2), (4),
(5),(6), \, \textrm{and} \, (10)$ may be found in~\cite{Taylor}.  We
refer the reader to these texts for original references and the
history of these results.  Specifically, the properties threshold and
a total vicinal preorder are two of eight equivalent conditions
presented in~\cite[Theorem 1.2.4]{Mahadev}.  The equivalence of
threshold graphs with unique realizability and degree-maximality
appears in~\cite[\S3.2]{Mahadev} along with six other conditions
determining degree sequences of threshold graphs.  The polytope of
graphical degree sequences is discussed in~\cite[\S3.3]{Mahadev}.  The
results of~\cite{Taylor} are not limited to the $k=2$ case as
outlined below.

For $k \geq 3$, the equivalence of 
\begin{enumerate} 
\item[$\bullet$]
threshold families and vertices of
$D_n(k)$ appears as \cite[Theorem 2.5]{Murthy},
\item[$\bullet$]
threshold families and $CC_t$ appears as \cite[Theorem 2.4.2]{Taylor},
\item[$\bullet$]
unique realizability and $DCC_t$ appears as \cite[Theorem 5.2.5]{Taylor}, 
\item[$\bullet$]
shiftedness and the $RRST$ condition appears as \cite[Theorem 2.5]{Reiterman}, and
\item[$\bullet$]
shiftedness and having a total vicinal preorder appears as \cite[Theorem 1]{Klivans},
\end{enumerate}
while the implications
\begin{enumerate}
\item[$\bullet$]
threshold implies uniquely realizable appears as \cite[Corollary 2.6]{Murthy},
\item[$\bullet$]
threshold implies shifted is an old observation, e.g. \cite[\S3.3,3.4]{Taylor} or \cite[\S10]{Golumbic}, and
\item[$\bullet$]
degree-maximal implies shifted appears as \cite[Proposition 9.3]{Duval}.
\end{enumerate}
\end{remark}

\begin{proof}(of Theorem \ref{omnibus-theorem}) 

\vskip.2in
\noindent
{\bf Equivalences:}

For the proof of equivalence of \eqref{threshold}, \eqref{vertex}, \eqref{CC},
consider the vector configuration $\VV:=\{\chi_S\}_{S \in \binom{[n]}{k}}$.
Note that all of the vectors in $\VV$ lie on the affine hyperplane $h(x)=k$,
where $h$ is the functional in $(\R^n)^*$ defined by $h(x):=\sum_{i=1}^n x_i$,
that is, they form an {\it acyclic vector configuration}, corresponding to
an affine point configuration in the above affine hyperplane.
The theory of zonotopes \cite[\S9]{Eckhoff} tells us that for a subset $K \subset \VV$ of an
acyclic configuration of vectors, the following three conditions are equivalent:
\begin{enumerate}
\item[(1)] There exists a linear functional $w$ with
$w(v) >0$ for $v \in K$ and $w(v) < 0$ for $v \in \VV \setminus K$.
\item[(2)] The sum $\sum_{v \in K} v$ is a vertex of the zonotope generated by the
line segments $\{[0,v]\}_{v \in \VV}$.
\item[(3)] The decomposition $\VV= K \sqcup \left( \VV \setminus K \right)$ is a 
{\it non-Radon partition} in the sense that the two cones positively generated by $K$ and 
by $\VV \setminus K$ intersect only in the zero vector.
\end{enumerate}
It should be clear that the first two of these three conditions correspond
to a $k$-family $K$ being threshold or being a vertex of $D_n(k)$.  The cancellation
condition $CC_t$ for all $t$ corresponds to the non-Radon partition condition
as follows.  The partition is non-Radon if one cannot have a dependence
$
\sum_{v \in K} a_v v = \sum_{v \in \VV \setminus K} b_v v 
$
with positive reals $a_v, b_v$.  Since $\VV=\{\chi_S\}_{S \in \binom{[n]}{k}}$
contains only integer vectors, without loss of generality the coefficients $a_v, b_v$ in such a
dependence can be assumed rational, and then by clearing denominators,
they may be assumed to be (positive) integers.  Furthermore, since $h(v)=k$ for each $v \in \VV$,
one may assume it is a {\it homogeneous} dependence, that is,
$$
\sum_{v \in K} a_v = \sum_{v \in \VV \setminus K} b_v \left( =:t \right).
$$
This homogeneous dependence now corresponds to a pair of $t$-tuples $(A_1,\ldots,A_t)$, $(B_1,\ldots,B_t)$
in which $a_v, b_v$ are the multiplicities of the sets, contradicting condition $CC_t$.

For the proof of equivalence of \eqref{uniquely-realizable} and \eqref{DCC},
note that if two $k$-families $K \neq K'$
had $d(K)=d(K')$ then the sets $\{A_1,\ldots,A_t\}:=K \setminus (K \cap K')$
and $\{B_1,\ldots,B_t\}:=K' \setminus (K \cap K')$ would contradict $CC_t$.  Conversely,
if one had two collections of $t$ sets $\{A_1,\ldots,A_t\} \subset K$ and
$\{B_1,\ldots,B_t\} \subset \VV \setminus K$ with a a dependence 
$\sum_{i=1}^t \chi_{A_i} = \sum_{i=1}^t \chi_{B_i}$, then 
$$
K':=\left( K \setminus \{A_1,\ldots,A_t\} \right) \cup \{B_1,\ldots,B_t\}
$$
would have $d(K') = d(K)$ but $K' \neq K$.

For the equivalence of \eqref{vicinal}, \eqref{RRST}, and \eqref{shifted}, one can easily check that
for a shifted family $K$ on $[n]$ one has $N_K[i] \supseteq N_K(j)$ whenever $i < j$, so
that the vicinal preorder $\vicin_K$ is total.  Conversely, if $\vicin_K$ is total,
relabel the set $[n]$ so that the integer order $<_\ZZ$ on $[n]$ is consistent with
$\vicin_K$ (that is, $i <_\ZZ j$ implies $i \vicin_K j$), and one can then check that
this labelling makes $K$ a shifted family.  Lack of totality for the vicinal preorder
means one has a pair $i,j$ with $i \not\vicin_K j$ 
(witnessed by some $(k-1)$-set $A$ in $N_K(j) \setminus N_K[i]$) 
and $j \not\vicin_K i$ 
(witnessed by some $(k-1)$-set $B$ in $N_K(i) \setminus N_K[j]$).  One can check that this
is exactly the same as a tuple $(A,B,i,j)$ which witnesses failure of the $RRST$ condition.

\vskip.2in
\noindent
{\bf Forward implications:}

To see \eqref{positive-threshold} implies \eqref{threshold}, note that for any $k$-set $S$,
the vector $x=\chi_S$ satisfies the inhomogeneous inequality $\sum_{i=1}^n c_i x_i > t$ if and only if
it satisfies the homogeneous inequality $w(x):=\sum_{i=1}^n \left(c_i - \frac{t}{k}\right) x_i> 0$.

It should also be clear that \eqref{CC} implies \eqref{DCC}, that is, $CC_t$ implies $DCC_t$.  
However, we emphasize the geometric statement\footnote{and see also 
\cite[Corollary 2.6]{Murthy} for a similar argument via linear programming.}
underlying the implication \eqref{threshold} implies \eqref{uniquely-realizable}:
when a subset $K$ of an acyclic 
vector configuration $\VV$ sums to a vertex of the zonotope generated by $\{[0,v]\}_{v \in \VV}$, no
other subset $K' \subset \VV$ can sum to the same vertex, otherwise the
dependence $\sum_{ v \in K \setminus (K \cap K')} v = \sum_{ v \in K' \setminus (K \cap K')} v$
would contradict the fact that $\VV = K \sqcup (\VV \setminus K)$ is a non-Radon partition
(or it would contradict the existence of the 
functional $w$ which is positive on $K$ and negative on $\VV \setminus K$).

To show \eqref{degree-maximal} implies \eqref{shifted}, we show the contrapositive.
Assume that $K$ is not isomorphic to any shifted family, and
without loss of generality, relabel $[n]$ so that $d(K)$ is weakly decreasing.
Since $K$ is still not shifted, by Proposition~\ref{prop:swing}(ii) it has
 applicable some swing that produces a family $K'$, for which $d(K')$ strictly majorizes $d(K)$.

To show \eqref{uniquely-realizable} implies \eqref{degree-maximal}, note that
\eqref{uniquely-realizable} is equivalent to \eqref{DCC} which implies $DCC_2$ and hence
\eqref{RRST}, which is equivalent to \eqref{shifted}.  Hence if $K$ is uniquely realizable, 
it is isomorphic to a shifted family.  Relabel so that $K$ is itself shifted, and hence
$d(K)$ is weakly decreasing.  
Choose a degree-maximal family $K'$ with $d(K') \unrhd d(K)$.  Then Proposition~\ref{prop:swing}(iii)
implies there exists  a $k$-family $K''$ with degree sequence $d(K'')=d(K)$ such that $K'$ can be
obtained from $K''$ by a (possibly empty) sequence of swings.  Unique realizability forces $K''=K$.
Since $K$ is shifted, there are no swings applicable to it (Proposition~\ref{prop:swing}(i)),
so the aforementioned sequence of swings must be empty, i.e. $K=K'$.  Hence $K$ is degree-maximal.

\vskip.2in
\noindent
{\bf Completing the circle of equivalences for $k=2$.}\
To show \eqref{shifted} implies \eqref{positive-threshold} when $k=2$, it
suffices to exhibit for any shifted $2$-family (graph) $K$ on vertex set $[n]$,
some positive coefficients $a_1,\ldots,a_n$ and threshold value $t$
such that $\{i,j\} \in K$  if and only if $a_i + a_j > t$.
Chv\'atal and Hammer \cite[Fact 3]{Chvatal} prove this by induction on $n$, as follows.  

  It is easy to see that in a shifted graph $K$, there always exists a vertex $v$ 
whose deletion $K \backslash v$ is a shifted graph, and for which either $v$ is {\it isolated} in the sense 
that $\{i,v\} \not\in K$ for all $i$, or $v$ is a {\it cone vertex} in the
sense that $\{i,v\} \in K$ for all $i \neq v$.  
Assume by induction that there are coefficients $\hat{a}_i$ and a threshold value
$\hat{t}$ exhibiting $K \backslash v$ as a positive threshold graph.  
If $v$ is a cone vertex, then taking
$a_v = t=\hat{t}$ and $a_i=\hat{a}_i$ for $i \neq v$ exhibits $K$ as a positive
threshold graph.  If $v$ is an
isolated vertex, then without loss of generality one can first perturb the $\hat{a}_i$ and $\hat{t}$ to
be rationals, and then clear denominators to make them positive integers.  After this,
taking $t=2\hat{t}+1, a_v=1$ and $a_i=2\hat{a}_i$ for $i \neq v$ exhibits $K$ as a positive
threshold graph.

\vskip.2in
\noindent
{\bf Strictness of the implications for $k \geq 3$.}

To show \eqref{shifted} $\nRightarrow$ \eqref{degree-maximal}, 
Example 9.4 of~\cite{Duval} gave a family of examples starting at $k=3$ and $n=10$.  
We give here an example with $n=9$.  Let $K$ be the shifted $3$-family of $[9]$
  generated by $\langle 178, 239, 456 \rangle.$  This shifted family has
  degree sequence $d(K) = (23,16,16,12,12,12,7,7,3),$ and is not degree-maximal:
  the family $K'$ generated by
  $\langle 149, 168, 238, 257, 356 \rangle$ has degree sequence
  $(23,17,15,13,12,11,8,6,3) \rhd d(K)$.

To show \eqref{degree-maximal} $\nRightarrow$ \eqref{uniquely-realizable}, check
the two shifted families   
$$
\begin{aligned}
&\langle 457, 168, 149, 248, 239 \rangle\\
&\langle 456, 357, 348, 267, 159 \rangle
\end{aligned}
$$
both achieve the degree-maximal degree sequence $d = (23,19,18,17,15,12,11,7,4)$. 

To show \eqref{uniquely-realizable} $\nRightarrow$ \eqref{threshold},
in terms of cancellation conditions, this amounts to showing that
satisfying $DCC_t$ for all $t$ is not equivalent to satisfying
$CC_t$ for all $t$, which is illustrated by an example in Theorem 5.3.1
of~\cite{Taylor}.

To show \eqref{threshold} $\nRightarrow$ \eqref{positive-threshold},
 in \cite[Example 2.1]{Reiterman} the authors observe that the $k$-family
 $H_m^k$ of $k$-subsets of $V:=\{-m,-m+1,\ldots-1,0,1,\ldots,m-1,m\}$
 given by 
 $$
 H_m^k:=\left\{ S \in \binom{V}{k}: \sum_{ s \in S } s > 0 \right\}
 $$
 is threshold for every $m, k$:  the functional $w(x)=\sum_{i=-m}^{m} i x_i$
 in $(\R^{|V|})^*$ has $w(\chi_S) >0$ for some subset $S \in \binom{V}{k}$
 if and only if $S \in H_m^k$, by definition.  They then show that $H_m^k$ is {\it not}
 positive threshold for $k\geq 3$ and $m$ larger than some bound $m_k$
 (with $m_3=7$).  A direct  example to show \eqref{shifted} $\nRightarrow$
 \eqref{threshold} is also given in \cite[Example 2.2]{Reiterman}.

\end{proof}

\begin{remark}\rm \
We have seen that shifted families do not always have uniquely realizable
degree sequences.  However, one might wonder whether it is possible for
a shifted family $K$ and a non-shifted family $K'$ to have the same
degree sequence.  This can happen, and follows from the method used to
prove \eqref{uniquely-realizable} implies \eqref{degree-maximal}, as we explain here.

Begin with a shifted family $K$ which is not degree-maximal,
and choose a degree-maximal family $K''$ with $d(K'') \unrhd d(K)$.
Then use Proposition~\ref{prop:swing}(iii) to find a family $K'$ with $d(K')=d(K)$
which has applicable a sequence of swings bringing it to $K''$.  Since $K$ is shifted
and therefore has no applicable swings, one knows $K' \neq K$.
\end{remark}

\begin{remark} \rm \
Other concepts for $k$-families have been considered, which 
may lie between threshold and shifted,
such as the $k$-families occurring in initial segments of {\it linear qualitative
probabilities} studied by Edelman and Fishburn \cite{EdelmanFishburn},
or the $k$-families inside {\it weakly} or 
{\it strongly acyclic games} \cite[Chapter 4]{Taylor}.

\end{remark}

\section{How to characterize degree sequences?}
\label{characterization-section}

\subsection{The problem, and an unsatisfactory answer}\

 Given a vector $d \in \PP^n$ with 
$$
|d|:=\sum_{i=1}^n d \equiv 0 \mod k,
$$
when is $d$ the degree sequence of some $k$-family $K$ on $[n]$?  

For $k=2$, there are many intrinsic characterizations of such {\it graphic sequences}; 
for example, see \cite{Sierksma} for $7$ such characterizations.   
The situation for $k$-families with $k > 2$ is
much less clear, although at least one has the following.

\begin{proposition}
\label{unsatisfactory-proposition}
The following are equivalent for a sequence $d = (d_1,\ldots,d_n) \in \PP^n$:
\begin{enumerate}
\item[(i)] $d=d(K)$ for some $k$-family.
\item[(ii)] There exists a degree-maximal $k$-family $K$ with $d(K) \unrhd d$.
\item[(iii)] There exists a shifted $k$-family $K$ with $d(K) \unrhd d$.
\end{enumerate}
\end{proposition}
\begin{proof}
As noted earlier, $d$ is majorized by its weakly decreasing rearrangement,
and so one may assume (by relabelling) that $d$ is weakly decreasing.

Then (i) implies (ii) trivially, (ii) implies (iii) because degree-maximality 
implies shiftedness, and (iii) implies (i) via 
Proposition~\ref{prop:swing}(iii).
\end{proof}

Unfortunately, the above proposition is an unsatisfactory answer, partly
due to the lack of an intrinsic characterization of degree sequences for shifted families,
or for degree-maximal families.

\begin{openproblem}
\label{characterization-problem}
  For $k \geq 3$, find simple intrinsic characterizations of degree sequences of
  \begin{enumerate}
   \item[(i)] $k$-families.
   \item[(ii)] degree-maximal $k$-families.
   \item[(iii)] shifted $k$-families.
  \end{enumerate}
\end{openproblem}

\begin{remark} (``Holes'' in the polytope of degree sequences?) \rm \ 

Fixing $k$ and $n$, there is an obvious inclusion
\begin{equation}
\label{obvious-inclusion}
\left\{\text{degrees }d(K)\text{ of }k\text{-families on }[n]\right\}
\,\, \subseteq \,\,
\left\{ d \in \N^n \cap D_n(k): \sum_{i=1}^n d_i \equiv 0\mod k \right\}
\end{equation}
where one should view the set on the right as the relevant lattice points
inside the polytope $D_n(k)$ which is the convex hull of degree sequences
of $k$-families.

For $k=2$, Koren \cite{Koren} showed that this inclusion is an equality.
He showed that the Erd\"os-Gallai {\it non-linear} inequalities, 
$$
\sum_{i=1}^k d_i \leq k(k-1) + \sum_{i=k+1}^n \min \{k,d_i\} \text{ for }k=1,2,\ldots,n,
$$
which give one characterization of degree sequences among all sequences 
$d$ of positive integers with even sum, 
are equivalent to
the following system of {\it linear} inequalities: 
$$
\sum_{i \in S} d_i - \sum_{j \in T} d_j \leq |S|(n-1-|T|)
$$
as $S,T$ range over all pairs of disjoint subsets $S,T \subseteq [n]$ 
which are not both empty.

Bhanu Murthy and Srinivasan \cite{Murthy} study several properties of
the polytope $D_n(k)$ for $k \geq 2$, including a description of some
of its facets, its $1$-skeleton, and its diameter.  

  There has been some speculation that for
$k \geq 3$ the inclusion \eqref{obvious-inclusion} is proper, that
is, there are non-degree sequence ``holes'' among the $d \N^n$ lattice points 
that lie within the convex hull $D_n(k)$ of degree sequences.  
However, we know of no such example, and have been able to check\footnote{by
comparing the values in Table~\ref{vertex-degree-sequence-table} with
values from a computer implementation of Stanley's method of lattice
point enumeration within a zonotope from \cite[\S 3]{Stanley-zonotope}.} 
that no such holes are present for $k=3$ and $n \leq 8$.

\begin{openproblem}
Are there ``holes'' in the polytope $D_n(k)$ of $k$-family vertex-degree
sequences?
\end{openproblem}

\end{remark}

\subsection{Some data on degree sequences}

Table~\ref{vertex-degree-sequence-table} lists some known
data on the number of vertex-degree sequences $d(K)$ for $k$-families $K$ on $[n]$, 
compiled via three sources:
\begin{enumerate}
\item[$\bullet$] The trivial values where $n \leq k+1$, along with
the equality of values for $(k,n)$ and $(n-k,n)$ that follows from the
second symmetry in Section~\ref{symmetries-subsection} below.
\item[$\bullet$] The values for $k=2$, which can be computed for
large $n$ explicitly using Stanley's results \cite{Stanley-zonotope}
on graphical degree sequences.
\item[$\bullet$] Brute force computation for $k=3$ by finding all
shifted $3$-families on $[n]$, computing all partitions majorized
by their degree sequences, and then summing the number of rearrangements
of these degree partitions.
\end{enumerate}

%%%%%%%%%%% TABLE  %%%%%%%%%%%%%%
% \begin{table}[hptb]
% \begin{center}
\begin{equation}
\label{vertex-degree-sequence-table}
\begin{tabular}{| l| l | l | l | l | l | l | l | l | l | l | l |} \hline
        &n=0  &n=1 &n=2&n=3&n=4&n=5&n=6  &n=7     &n=8\\ \hline 
 k = 1  &1    & 2  & 4 & 8 & 16& 32& 64  & 128    &256\\ \hline
 k = 2  &1    & 1  & 2 & 8 & 54&533&6944 & 111850 &2135740\\ \hline
 k = 3  &1    & 1  & 1 & 2 & 16&533&42175& 5554128&1030941214\\ \hline
 k = 4  &1    & 1  & 1 & 1 & 2 & 32&6944 & 5554128&?   \\ \hline
 k = 5  &1    & 1  & 1 & 1 & 1 & 2 & 64  & 111850 &1030941214\\  \hline
 k = 6  &1    & 1  & 1 & 1 & 1 & 1 & 2   & 128    &2135740\\  \hline
 k = 7  &1    & 1  & 1 & 1 & 1 & 1 & 1   & 2      &256 \\ \hline
 k = 8  &1    & 1  & 1 & 1 & 1 & 1 & 1   & 1      &2   \\ \hline
\end{tabular}
\end{equation}
% \end{center}
% \end{table}
%%%%%%%%%%%%%%%%%%%%%%%%%%%%%%%%%%%%%%%%%%%%%%%%%%

\subsection{Reconstructing families}

  In preparation for what follows, we note some easy facts
about reconstructing $k$-families and shifted $k$-families from other data.

\begin{definition}
Extend the open neighborhood notation $N_K(i)$ from Definition~\ref{vicinal-definition}
as follows.  Given a $k$-family $K$ on $[n]$ and any subset $T$ of $[n]$, the 
{\it open neighborhood} or {\it link} of $T$ in $K$ is
$$
N_K(T):=\{S \subset [n]: S \cup T  \in K \text{ and } S \cap T = \emptyset \}.
$$
Note that if $|T|=i$ then for $S \in N_K(T)$, $|S| = k-i$.  Also note that when $|T|=k$, one either
has 
$$
N_K(T) = \begin{cases} \{\emptyset\} & \text{ if }T \in K \\
                         \emptyset\  & \text{ if }T \not\in K.
         \end{cases}
$$
Define the {\it $i$-degree sequence/function} $d_K^{(i)}$ as follows:
$$
\begin{aligned}
d_K^{(i)} : & \binom{[n]}{i}  &\longrightarrow  & \quad \N \\
            &    T           & \longmapsto     & \quad |N_K(T)|.       
\end{aligned}
$$
Thus the vertex-degree sequence $d(K)=(d_1(K),\ldots,d_n(K))$ 
is essentially the same as the function $d^{(1)}_K: [n] \rightarrow \N$.
\end{definition}

\begin{proposition}
\label{reconstruction-trivialities}
Let $K$ be a $k$-family on $[n]$.
\begin{enumerate}
\item[(i)] For any $i$ in the range $0 \leq i \leq k$, the restricted function 
$N_K(-): \binom{[n]}{i} \rightarrow \binom{[n]}{k-i}$ determines $K$ uniquely.
\item[(ii)] If $K$ is a shifted $k$-family, and $T$ any $i$-subset of $[n]$, then $N_K(T)$ is a
shifted $(k-i)$-family on the set $[n] \setminus T$ (linearly ordered in the usual way).
\item[(iii)] If $K$ is a shifted $k$-family, then the subfacet-degree function 
$d_K^{(k-1)}: \binom{[n]}{k-1} \rightarrow \N$ determines $K$ uniquely.
\end{enumerate}
\end{proposition}

\begin{proof}
Only assertion (iii) requires comment.  By assertion (i), one only needs to check that,
when $K$ is a shifted $k$-family, the (set-valued) function $N_K(-)$ restricted
to $(k-1)$-sets is determined uniquely by the (integer-valued) function $d_K^{(k-1)}$.  But assertion (ii)
implies that for any $(k-1)$-subset $T$, the collection $N_K(T)$ is a shifted $1$-family 
on $[n] \setminus T$, which implies that $N_K(T)$
is completely determined by its cardinality, namely the integer $d_K^{(k-1)}(T)$
\end{proof}

\begin{remark}
This proposition perhaps suggests that results about the vertex-degree
sequence $d(K)=d^{(1)}_k$ for $k=2$ could generalize in different directions:
in the direction of vertex-degree sequences $d_K^{(1)}$, or in
the direction of {\it subfacet} (=$(k-1)$-set) degree sequences $d_K^{(k-1)}$.
\end{remark}

\begin{remark}
We wish to underscore a difference between vertex-degrees and subfacet-degrees. 

It is natural to identify vertex-degree functions $d^{(1)}_K: [n] \rightarrow \N$ with 
the vertex-sequences $d(K)=(d_1(K),\ldots,d_n(K))$.  Furthermore,
it suffices to characterize those which are weakly decreasing;  this just
means characterizing the degree functions up to the natural
action of the symmetric group $\Sym_n$ on the domain $[n]$ of $d^{(1)}_K$.

For subfacet-degree functions $d^{(k-1)}_K$, however, it is not so natural to
identify them with some sequence of degree values, as this involves the choice
of a linear ordering on the domain $\binom{[n]}{k-1}$ to write down such a sequence.
But it is still true that it suffices to characterize 
subfacet-degree functions $d^{(k-1)}_K$ up to the action 
of $\Sym_n$ on their domain $\binom{[n]}{k-1}$.
\end{remark}

\subsection{Some promising geometry}

  The goal of this subsection is to shed light on 
Open Problem~\ref{characterization-problem}(iii),
motivated by a characterization of graphical degree sequences (i.e. $k=2$) 
due to Ruch and Gutman \cite{Ruch}, and reformulated
by Merris and Roby \cite{MerrisR}.

Given a weakly decreasing sequence $d \in \PP^n$, identify $d$ with its {\it Ferrers
diagram} as a partition, that is the subset of boxes 
$
\{(i,j) \in \PP^2: 1 \leq j \leq d_i, \,  1 \leq i \leq n \}
$
in the plane $\PP^2$.  The {\it conjugate} or {\it transpose} partition $d^T$ 
is the one whose Ferrers diagram is obtained by swapping $(i,j)$ for $(j,i)$,
and the {\em trace} or {\em Durfee rank} of $d$ is the number of diagonal boxes
of the form $(i,i)$ in its Ferrers diagram, that is, 
$\mathrm{trace}(d) = |\{ j : d_j \geq j \}|$.  Ruch and Gutman's characterization
says the following.

\begin{theorem}[\cite{Ruch}]
\label{thm:RG}
An integer sequence $d \vdash 2m $ is the degree sequence of some $2$-family if and only if
$$ 
\sum_{i=1}^{k} (d_i + 1) \leq \sum_{i=1}^{k} {d_i}^T, \, \, \, \, \,
1 \leq k \leq \textrm{\em trace}(d). 
$$
\end{theorem}

\begin{figure}[ht]
\centering
\includegraphics[height=1.6in]{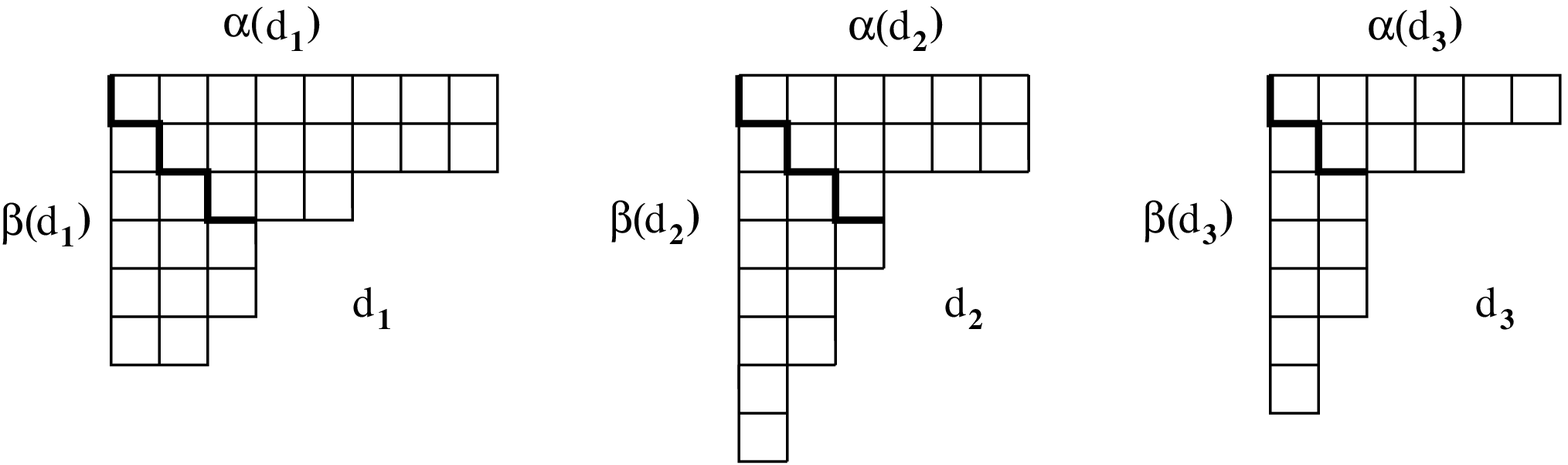}
\caption{Cutting partitions into shifted shapes.}
\label{fig:Fer}
\end{figure}

\noindent
Merris and Roby's reformulation of this result uses some geometry
of diagrams for strict partitions placed in the {\it shifted plane}, which is
the set of boxes lying weakly above the diagonal in the usual 
positive integer plane $\PP^2$. Given a Ferrers diagram $d$ embedded in the usual
plane, cut it into two pieces along the
``subdiagonal staircase'' as shown in Figure~\ref{fig:Fer}.  
Let $\alpha(d)$ (resp. $\beta(d)$) denote the subshape formed by the boxes with $i\leq j$ (resp. $i > j$).
If the trace of $d$ is $t$ then the sequence of row sizes of $\alpha(d)$ form a
{\it strict partition} $\alpha_1 > \cdots > \alpha_t$ with $t$ parts, that we will also denote
by $\alpha(d)$.  Similarly the column sizes of $\beta(d)$ form a strict partition with $t$ parts
that we will denote $\beta(d)$.

Here is the Merris and Roby formulation. Recall that $\beta \succ \alpha$ means
that $\beta$ weakly majorizes $\alpha$.

\begin{theorem} (Theorem 3.1~\cite{MerrisR})
\label{thm:MR}
An integer sequence $d \vdash 2m$ is graphical, that is, $d(K)$ for some
$2$-family $K$, if and only if $\beta(d) \succ \alpha(d) $.  

Moreover, $d=d(K)$ for a shifted family $K$ if and only if the strict partitions
$\alpha(d) , \beta(d)$ are the same.  
\end{theorem}

\begin{example}
In Figure~\ref{fig:Fer}, ${\beta}(d_1) \nsucc {\alpha}(d_1) $, ${\beta}(d_2)
\triangleright {\alpha}(d_2)$, and ${\beta}(d_3) = {\alpha}(d_3)$.  Therefore $d_1$ is not
graphical, $d_2$ is graphical but not shifted, and $d_3$ is
shifted.
\end{example}

We codify here some of the geometry relating the plane $\PP^2$ and 
the shifted plane $\binom{\PP}{2}$
that makes this work, and which will generalize in two directions to $k$-families for
arbitrary $k$.  Given integers $i < j$, let $[i,j]:=\{i,i+1,\ldots,j-1,j\}$.
Consider the following decomposition of the finite rectangle $[1,n-1] \times [1,n]$
inside $\PP^2$:
$$
[1,n-1] \times [1,n] : = \binom{[n]}{2} \sqcup f \binom{[n]}{2}
$$
where we are identifying $\binom{[n]}{2}$ with
$\{(i_1,i_2) \in \N^2: 1 \leq i_1 < i_2 \leq n\}$
and where $f: \R^2 \rightarrow \R^2$ is the affine isomorphism $f(i_1,i_2) := (i_2,i_1)+(-1,0)$.

Given a $2$-family/graph $K$ on $[n]$, think of $K$ as a subset of $\binom{[n]}{2}$,
and let 
$$
\pi(K):= K \sqcup f(K)  \subset [1,n-1] \times [1,n].
$$
We further rephrase the notation of Theorem~\ref{thm:MR}.
Relabel so that $d(K)$ is weakly decreasing, and consider the (French-style)
Ferrers diagram $\pi \subset [1,n-1] \times [1,n]$
which is left-justified and has $d_i(K)$ cells with $x_2$-coordinate $i$.
Let $\alpha(K):=\pi \cap   \binom{[n]}{2}$ 
and $\beta(K):= \pi \cap f \binom{[n]}{2} $.

\begin{proposition}
\label{Merris-Roby-essence}
For any $2$-family $K$ on $[n]$ one has the following.
\begin{enumerate}
\item[(i)] The degree sequence $d(K)=(d_1(K),\ldots,d_n(K))$ 
has $d_i(K)$ given by the number of boxes in $\pi(K)$ with $x_2$-coordinate equal to $i$,
and this sequence completely determines $K$ if $K$ is a shifted $2$-family.
\item[(ii)] The following are equivalent 
 \begin{enumerate}
  \item[(a)] $K$ is shifted, that is, it forms a componentwise order 
ideal in $\binom{[n]}{2}$.
  \item[(b)] $\pi(K)$ is a componentwise order ideal of $[1,n-1] \times [1,n]$.
  \item[(c)] $\alpha(K) = f^{-1}(\beta(K))$, and both coincide with
the family $K$, thought of as a subset of $\binom{[n]}{2}$ inside $\PP^2$.
 \end{enumerate}
\end{enumerate}
\end{proposition}

To generalize some of this for arbitrary $k$, we recall a well-known triangulation 
of the prism over the $(k-1)$-simplex $\{ x \in \R^{k-1}: 0 \leq x_1 \leq \cdots \leq x_{k-1} \leq 1\}$:
$$
\{ x \in \R^k: 0 \leq x_1 \leq \cdots \leq x_{k-1} \leq 1\text{ and }0 \leq x_k \leq 1\} \\
 = \sigma_1 \cup \cdots \cup \sigma_k
$$
where 
$$
\sigma_k=\sigma:=\{ x \in \R^k:0 \leq x_1 \leq \cdots \leq x_{k} \leq 1 \}
$$
is a $k$-simplex, and for $j=1,2,\ldots,k-1$,
$$
\begin{aligned}
\sigma_j &:= \{ x \in \R^k:
              0 \leq x_1 \leq \cdots \leq x_{j-1} \leq x_k 
                  \leq x_j \leq x_{j+1} \leq \cdots \leq x_{k-1} \leq 1\}\\
         &= f_j(\sigma)
\end{aligned}
$$
where $f_j: \R^k \rightarrow \R^k$ is the linear isomorphism 
$$
f_j(i_1,\ldots,i_k) := (i_1,i_2,\ldots,i_{j-1},i_{j+1},i_{j+2},\ldots,i_{k-1},i_k,i_j).
$$

\begin{figure}[ht]
\centering
\includegraphics[height=2in]{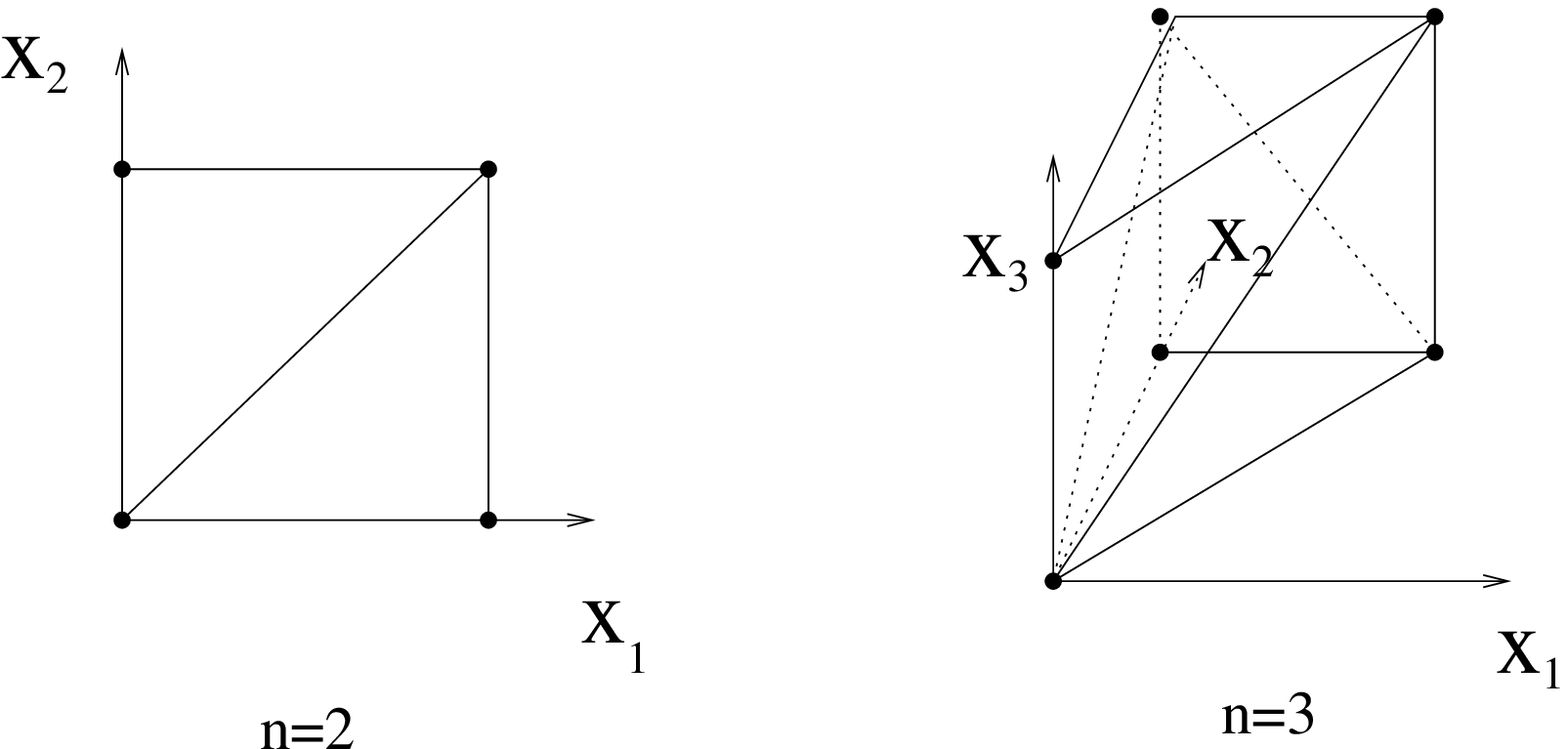}
\caption{A well-known triangulation of a prism over a simplex.}
\label{prism-triangulation}
\end{figure}

This triangulation is depicted for $n=2,3$ in Figure~\ref{prism-triangulation}.
It arises, for example,
\begin{enumerate}
\item[$\bullet$]
as a convenience in proving facts about homotopies in simplicial sets \cite{May},
\item[$\bullet$]
as the special case of the {\it staircase triangulation} of a product of simplices
\cite[Example 8.12]{Sturmfels},
where one of the simplices is $1$-dimensional, or
\item[$\bullet$]
as the special case of the {\it $P$-partition triangulation} of the {\it order polytope} 
\cite{Stanley-order-polytope}, where the poset $P$ is the disjoint union of 
two chains having sizes $k-1$ and $1$.
\end{enumerate}

We wish to apply this triangulation toward
understanding vertex-degree functions $d_K^{(1)}=d(K)$ and subfacet-degree functions $d_K^{(k-1)}$ of 
$k$-families and shifted $k$-families on $[n]$.  For this, we dilate the triangulation by $n$, 
and consider two different ways to decompose the lattice points within these (dilated) objects.

\begin{definition}
Fix $n$ and $k$, and identify
$$
\binom{[n]}{k} = \{ (i_1,\ldots,i_k) \in \PP^k : 1 \leq i_1 < \cdots < i_k \leq n\}.
$$

The {\it vertex-degree decomposition} of 
$$
\binom{[n-1]}{k-1} \times [1,n] 
   := \{ x \in \PP^k: 1 \leq x_1 < \cdots < x_{k-1} \leq n-1, 
             \text{ and } 1 \leq x_k \leq n \}
$$
is $\binom{[n-1]}{k-1} \times [1,n]  = \sigma^\vertex_1 \sqcup \cdots \sqcup \sigma^\vertex_n$,
in which 
$$
\begin{aligned}
\sigma^\vertex_k &:=\binom{[n]}{k} \\
\sigma^\vertex_j &:= \{ 1 \leq x_1 <\cdots < x_{j-1} < x_k  \leq x_j < x_{j+1} < \cdots < x_{k-1} \leq n-1\} \\
       &= f^\vertex_j\binom{[n]}{k} \text{  for }j=1,2,\ldots,k-1
\end{aligned}
$$
where $f^\vertex_j(x) :=f_j(x) + (\underbrace{0,\ldots,0}_{j-1\text{ positions}},
                               \underbrace{-1,\ldots,-1}_{k-j\text{ positions}},0).$

The {\it subfacet-degree decomposition} of  
$$
\binom{[n]}{k-1} \times [k,n] 
    := \{ x \in \PP^k: 1 \leq x_1 < \cdots < x_{k-1} \leq n, 
             \text{ and } k \leq x_k \leq n \}
$$
is
$\binom{[n]}{k-1} \times [k,n]  = \sigma^\subfacet_1 \sqcup \cdots \sqcup \sigma^\subfacet_n$,
in which 
$$
\begin{aligned}
\sigma^\subfacet_k &:=\binom{[n]}{k}\\
\sigma^\subfacet_j &:= \{ 1 \leq x_1 <\cdots < x_{j-1} 
                            < x_k-(k-j) < x_j < x_{j+1} < \cdots < x_{k-1} \leq n\} \\
       &= f^\subfacet_j\binom{[n]}{k}\text{ for }j=1,2,\ldots,k-1
\end{aligned}
$$
where $f^\subfacet_j(x) :=f_j(x) + (0,\ldots,0,k-j).$

Note that $f_k$ is the identity map, and 
the formulae for $ f^\vertex_j, f^\subfacet_j$ are consistent
with defining both $f^\vertex_k, f^\subfacet_k$ also as the identity map.
\end{definition}

We omit the straightforward verification of the following:
\begin{proposition}
The vertex-degree and subfacet-degree decompositions really are disjoint decompositions
of the claimed sets, $\binom{[n-1]}{k-1} \times [1,n]$ and $\binom{[n]}{k-1} \times [k,n]$.
\end{proposition}

\begin{definition}
For a $k$-family $K$ on $[n]$, thinking of $K$ as a subset of 
$\binom{[n]}{k}= \sigma^\vertex_k=\sigma^\subfacet_k$, define 
$$
\begin{aligned}
\pi^\vertex(K) & := \sqcup_{j=1}^k f^\vertex_j(K)  &\subset &\binom{[n-1]}{k-1} \times [1,n]\\
\pi^\subfacet(K)& := \sqcup_{j=1}^k f^\subfacet_j(K) & \subset &\binom{[n]}{k-1} \times [k,n].
\end{aligned}
$$
\end{definition}

\begin{figure}[ht]
\centering
\includegraphics[height=2in]{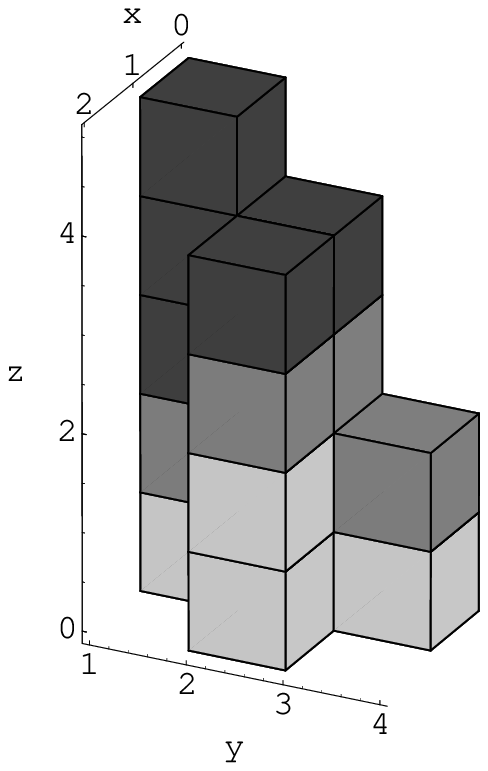} \hspace{.3in}
\includegraphics[height=1.6in]{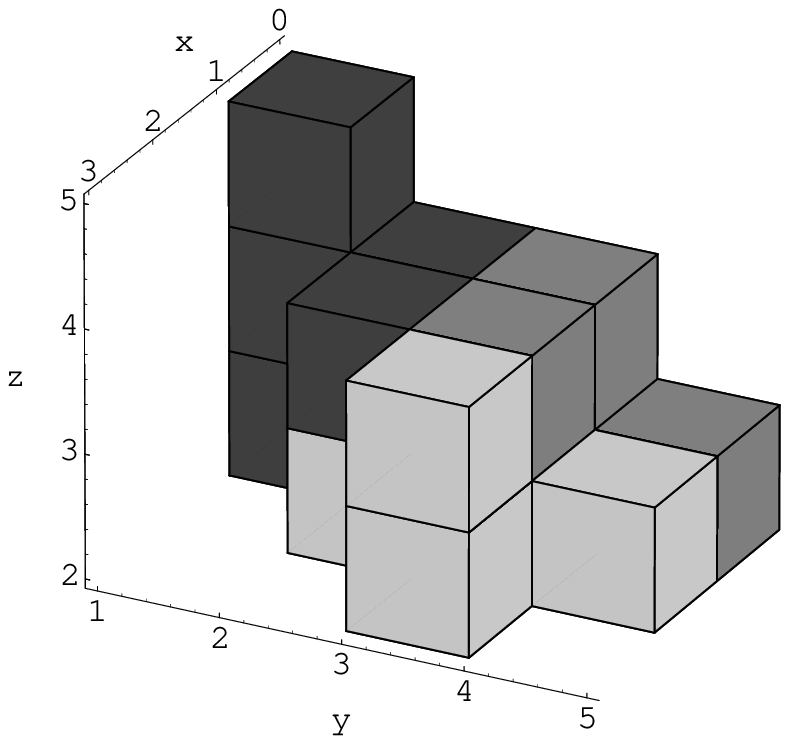}
\caption{$\pi^\vertex(K)$ and $\pi^\subfacet(K)$ from Example~\ref{ex:decomposition}}
\label{fig:decomposition}
\end{figure}

\begin{example}
\label{ex:decomposition}
Let $K$ be the shifted $3$-family, $K = \{123, 124, 134, 234, 125\}$.
Figure~\ref{fig:decomposition} shows $\pi^\vertex(K)$ and
$\pi^\subfacet(K)$.  $f_1^\vertex(K)$ and $f_1^\subfacet(K)$ are the collections
of lightest shaded cubes and $f_3^\vertex(K)$ and $f_3^\subfacet(K)$ are the collections of darkest
shaded cubes.

\end{example}

The key point of these constructions is their analogy to the $\alpha,
\beta$ appearing in Theorem~\ref{thm:MR} and
Proposition~\ref{Merris-Roby-essence}.
Proposition~\ref{k-family-geometry} below generalizes some of their
assertions, but requires a little more notation.

\begin{definition}
Let 
$$
\begin{aligned}
\rho_{k-1}&: \R^k &\rightarrow &\R^{k-1} \\
\rho_1&: \R^k &\rightarrow &\R^1
\end{aligned}
$$
denote the usual orthogonal projections from 
$\R^k$ onto its first $k-1$ coordinates and its last coordinate.

Given a subset $\pi$ of $\binom{[n]}{k-1} \times [k,n]$, let $\lambda(\pi)$ be the
unique subset of  $\binom{[n]}{k-1} \times [k,n]$ obtained by ``pushing the cells of $\pi$
down in the $x_k$ coordinate'', that is, for each $(k-1)$-set $T$, the
fiber intersection $\rho_{k-1}^{-1}(T) \cap \lambda(\pi)$ should have 
\begin{enumerate}
\item[$\bullet$] the same cardinality as the intersection  $\rho_{k-1}^{-1}(T) \cap \pi$, but
\item[$\bullet$] its $x_k$-coordinates forming an initial segment of $[k,n]$.
\end{enumerate}
\end{definition}

\begin{proposition} 
\label{k-family-geometry}
Let $K$ be any $k$-family on $[n]$.
\begin{enumerate}
\item[(i)] The subfacet-degree 
function $d_K^{(k-1)}:\binom{[n]}{k-1} \rightarrow \N$ is given by
$$
d_K^{(k-1)}(T) = | \rho_{k-1}^{-1}(T) \cap \pi^\subfacet(K) |.
$$
\item[(i$^\prime$)] The vertex-degree function $d_K^{(1)}:[n] \rightarrow \N$ is given by
$$
d_K^{(1)}(i) = | \rho_1^{-1}(i) \cap \pi^\vertex(K) |.
$$
\item[(ii)] Letting $\pi:=\pi^\subfacet(K)$ and $\lambda:=\lambda(\pi)$,
the following are equivalent 
 \begin{enumerate}
  \item[(a)] $K$ is shifted, that is, a componentwise order ideal of $\binom{[n]}{k}$.
  \item[(b)] $\pi^\subfacet(K)(=\pi)$ is a componentwise order ideal of $\binom{[n]}{k-1} \times [k,n]$.
  \item[(b$^\prime$)]  $\pi^\vertex(K)$ is a componentwise order ideal of $\binom{[n-1]}{k-1} \times [1,n]$.
  \item[(c)] The sets $ (f_j^\subfacet)^{-1}( \lambda \cap \sigma_j^\subfacet)$ are all equal to $K$ for 
      $j=1,2,\ldots,k$.
 \end{enumerate}
\end{enumerate}
\end{proposition}

\begin{proof}
\text{ }\\
\noindent
{\it Proof of (i) and (i$^\prime$)}.

Note that a set $S=\{i_1 < \cdots < i_k\}$ in $K$ has 
\begin{enumerate}
\item[$\bullet$]
the last coordinate of its image $f_j^\vertex(S)$ equal to $i_j$, while
\item[$\bullet$]
the first $k-1$ coordinates of its $f_j^\subfacet(S)$ equal to $(i_1,\cdots,\hat{i_j},\cdots,i_k)$.
\end{enumerate}
Thus the various images of $S$ under the maps  $f_j^\vertex, f_j^\subfacet$ contribute 
to the correct components of the appropriate degree sequences.

\noindent
{\it Proof of (ii)}.
First note that either of (b) or (b$^\prime$) implies (a):
since $K$ is the intersection of $\pi^\subfacet(K)$ (resp. $\pi^\vertex(K)$) with $\binom{\PP}{2}$,
this means $K$ will form an order ideal of $\binom{\PP}{2}$ if $\pi^\subfacet(K)$ (resp. $\pi^\vertex(K)$)
is a componentwise order ideal of  $\binom{[n]}{k-1} \times [k,n]$ (resp. $\binom{[n-1]}{k-1} \times [1,n]$).

To show (a) implies (b), assume $K$ is a shifted $k$-family on $[n]$, and
that we are given a vector $x=(x_1,\ldots,x_k) \in \pi^\subfacet(K)$.
We must show that if one lowers some coordinate of $x$ by one and the result $x'$ remains in
$\binom{[n]}{k-1} \times [k,n]$, then one still has $x' \in \pi^\subfacet(K)$.  
Let $j:=j(x)$ be the unique index in $1,2,\ldots,k$ such that $x \in \sigma^\subfacet_j$,
and say $x=f^\subfacet_j(i_1,\ldots,i_k)$ for $S=\{i_1,\ldots,i_k\} \in K$.

\noindent
{\bf Case 1.} $j(x')=j$.  

Then $x'=f_j^\subfacet(S')$ for some 
$k$-set $S' \in \binom{\PP}{k}$, and we wish to show that $S' \in K$,
so that $x' \in \pi^\subfacet(K)$.
Since $x'$ is componentwise below $x$, and the inverse 
map $(f_j^\subfacet)^{-1}$ is easily checked to be
componentwise order-preserving, $S'$ lies componentwise below 
$S$ in $\binom{\PP}{k}$.  Hence $S'$ is also in the shifted family $K$,
as desired.

\noindent
{\bf Case 2.} $j(x')\neq j$.  

By the definition of $j=j(x)$, one knows that 
$$
1 \leq x_1 < \cdots < x_{j-1} < x_k - (k-j) < x_j < \cdots < x_{k-1} \leq n.
$$
There are two possibilities for how $x'$ might fail to satisfy the same inequalities.

\noindent
{\bf Case 2a.} $x_k - (k-j) = x_j - 1= x'_j$.

In this case, one checks that $x'=f^\subfacet_{j+1}(S) \in \pi^\subfacet(K)$.

\noindent
{\bf Case 2b.} $x_{j-1} = x_k - (k-j) - 1 = x'_k - (k-j)$ (so $x'_k=x_k-1$).

In this case $x'=f^\subfacet_{j-1}(S')$
where $S'=\{i_1,\ldots,i_{j-2},i_{j-1}-1,i_j-1,i_{j+1},\ldots,i_k\}$.
Because $K$ is shifted, $S'$ is also in $K$, so $x' \in \pi^\subfacet(K)$.

The proof that (a) implies (b$^\prime$) is extremely similar.  Case 1 is the same, and
Case 2 breaks up into these two possibilities depending upon how $x'$ fails
to satisfy the inequalities satisfied by $x$
$$
1 \leq x_1 < \cdots < x_{j-1} < x_k \leq x_j < x_{j+1} \cdots < x_{k-1}
$$
that come from $j(x)=j$:

\noindent
{\bf Case 2a$^\prime$.} $x_{j-1} = x_k - 1=x'_k$.

In this case, $x'=f^\vertex_{j-1}(S) \in \pi^\vertex(K)$.

\noindent
{\bf Case 2b$^\prime$.} $x_k = x_j - 1 = x'_j - 2$ (so $x'_j=x_j-1$).

In this case $x'=f^\vertex_{j+1}(S')$
where $S'=\{i_1,\ldots,i_{j-1},i_{j}-1,i_{j+1}-1,i_{j+2},\ldots,i_k\}.$
Because $K$ is shifted, $S'$ is also in $K$, so $x' \in \pi^\vertex(K)$.

\vskip.1in
To show (b) implies (c), note that if $\pi^\subfacet(K)$ is a componentwise order ideal of 
$\binom{[n]}{k-1} \times [k,n]$, then $\lambda = \pi^\subfacet(K)$, that is, $\pi^\subfacet(K)$
has already been ``pushed down'' in the $x_k$ direction.  Thus
for every $j$ one has 
$$
\lambda \cap \sigma^\subfacet_j = \pi^\subfacet(K) \cap \sigma^\subfacet_j = f^\subfacet_j(K),
$$
and hence $\left(f^\subfacet_j\right)^{-1}(\lambda \cap \sigma^\subfacet_j)  = K$.

To show (c) implies (a), assume 
$\left(f^\subfacet_j\right)^{-1}(\lambda \cap \sigma^\subfacet_j)  =  K$ for each $j=1,2,\ldots,n$.
We wish to deduce $K$ is shifted, so given $S=\{i_1,\ldots,i_k\}$ with 
$1 \leq i_1 < \cdots < i_k \leq n$, it suffices to show that
if $S'=\{i_1,\ldots,i_{j-1},i_j-1,i_{j+1},\ldots,i_k\}$ still has $i_{j-1} < i_j - 1$ then
$S' \in K$.  Since $\left(f^\subfacet_j\right)^{-1}(\lambda \cap \sigma^\subfacet_j) = K$, one knows
$$
f^\subfacet_j(S)=(i_1,\ldots,\hat{i_j},\ldots,i_{k-1},i_j+(k-j)) \in \lambda.
$$
Hence
$$
(i_1,\ldots,\hat{i_j},\ldots,i_{k-1},i_j+(k-j)-1) \in \lambda
$$
since $\lambda$ is closed under lowering the $x_k$-coordinate, within the range $[k,n]$,
and $i_j+(k-j)-1$ is still in this range:
$$
i_j - 1 > i_{j-1} \geq j-1 \quad \Rightarrow \quad i_j+(k-j)-1 \geq k.
$$
But $(i_1,\ldots,\hat{i_j},\ldots,i_{k-1},i_j+(k-j)-1) = f^\subfacet_j(S')$,
so $S' \in \left(f^\subfacet_j\right)^{-1}(\lambda \cap \sigma^\subfacet_j) = K$.
\end{proof}

\begin{example}
The three non-shifted $3$-families
$$
\begin{aligned}
K_1 &=\{124\} \\
K_2 &=\{123,134\} \\
K_3 &=\{123,124,234\}
\end{aligned}
$$
illustrate the necessity of comparing {\it all} $k$ of the
sets $(f_j^\subfacet)^{-1}( \lambda \cap \sigma_j^\subfacet)$ in condition
(ii)(c) above.  For $K_1$, the sets for $j=1,2$ coincide with $K$,
but $j=3$ does not.  For $K_2$, the sets for $j=1,3$ coincide with $K$,
but $j=2$ does not.  For $K_3$, the sets for $j=2,3$ coincide with $K$,
but $j=1$ does not.

Note however, that all $3$ of these families $K_1,K_2,K_3$ are
isomorphic to shifted families, by reindexing the set $[n]=[4]$.
\end{example}

\begin{remark}
One might hope to
characterize $d^{(k-1)}_K$ for $k$-families $K$ by saying that
the sets $\alpha_1,\ldots,\alpha_k$ with
$$
\alpha_j:=(f_j^\subfacet)^{-1}(\lambda \cap \sigma_j^\subfacet)
$$
obey some inequalities with respect to some partial order generalizing the
weak  majorization $\alpha(d) \prec \beta(d)$ in
Theorem~\ref{thm:MR}.  Presumably, any such partial order would be stronger
than the ordering by weight (=number of cells), so we consider some examples to see
how the $\alpha$ can be ordered by weight.

The example $K_4=\{123,145\}$ has 
$$
\begin{aligned}
\alpha_1&=\{123,145\} \\
\alpha_2&=\{123,124,125\} \\
\alpha_3&=\{123\}
\end{aligned}
$$
so that $|\alpha_2|> |\alpha_1|, |\alpha_3|$.

On the other hand, $K_5=\{123,456\}$ has
$$
\begin{aligned}
\alpha_1&=\{123,145,146,156\}\\
\alpha_2&=\{123\} \\
\alpha_3&=\{123\}
\end{aligned}
$$
so that $|\alpha_1|> |\alpha_2|, |\alpha_3|$.

These would seem to preclude an assertion that the $\alpha_j$ are {\it totally}
ordered by something like a weak majorization.  One might be tempted to conjecture the following:

For any $k$-family, one has
$\alpha_k \prec \alpha_1, \alpha_2,\ldots,\alpha_{k-1}$,
where we generalize $\alpha \prec \beta$ to mean that for every order ideal $I$ of $\binom{[n]}{k-1}$
one has an inequality
\begin{equation}
\label{higher-majorization}
\sum_{S \in I} d^{(k-1)}(\alpha)(S) \leq 
\sum_{S \in I} d^{(k-1)}(\beta)(S).
\end{equation}

However, the family $K_6=\{123,124,135\}$ has
$$
\begin{aligned}
\alpha_1&=\{123, 124, 135\}\\
\alpha_2&=\{123, 125, 124\}\\
\alpha_3&=\{123, 124, 134\}
\end{aligned}
$$ 
and one can check $\alpha_3 \nprec \alpha_2$ because the order ideal $I$ generated by $\{34\}$ fails to satisfy the  inequality \eqref{higher-majorization}.

\end{remark}

\begin{remark}
For $k=3$, the associations between a shifted $3$-family $K \subset
\binom{[n]}{3}$ and the componentwise order ideals $\pi^\vertex(K),
\pi^\subfacet(K)$ in $\N^3$ are reminiscent of the correspondence used
in~\cite{Klivans} relating shifted families and {\it totally symmetric
plane partitions}.  In \cite{Klivans} this was used for the purposes
of enumerating shifted $3$-families.
\end{remark}

\section{Shifted families and plethysm of elementary symmetric functions}
\label{plethysm-section}

The goal of this section is to review the well-known equivalence between the 
study of degree sequences for $k$-families
and the problem of computing plethysms of the elementary symmetric functions,
as well as to push this a bit further.
We refer to the books by Macdonald \cite{Macdonald}, Sagan \cite{Sagan}, or
Stanley \cite[Chap. 7]{Stanley-EC2}
for symmetric function facts and terminology not defined here.

\begin{definition} \rm \  
Define a symmetric function $\Psi_k(\xx)$ in variables $\xx=(x_1,x_2,\ldots)$ 
by any of the first four equations below:

\begin{eqnarray}
\Psi_k(\xx)& :=& \sum_{k\text{-families }K} \xx^{d(K)} \\
          & :=& \prod_{S \in \binom{\PP}{k}} (1 + \xx_S) \\
          & :=& \sum_{\text{ partitions }\mu} c_{\mu,k} m_\mu \\ 
          & :=& \sum_{m \geq 0} e_m[e_k] \\
          & = &\sum_{\lambda: |\lambda| \equiv 0 \mod k} a_{\lambda,k} s_\lambda
\end{eqnarray}
\noindent
where here 
$$
\begin{aligned}
\xx^d &:=x_1^{d_1} x_2^{d_2} \cdots \\
\xx_S &:= \prod_{i \in S} x_i \\
c_{\mu,k} &:=|\{ k\text{-families }K\text{ realizing }d(K)=\mu \}| \\
m_\mu &:= m_\mu(\xx) = \text{ monomial symmetric function corresponding to }\mu \\
s_\lambda &:= s_\lambda(\xx) = \text{ Schur function corresponding to }\mu \\
e_m &:= \sum_{S \in \binom{\PP}{m}} \xx_S = m^{th} \text{ elementary symmetric function } \\
e_m[e_k] &:= \text{ plethysm of }e_k\text{ into }e_m \\
\end{aligned}
$$
and $a_{\lambda,k}$ are plethysm coefficients:  the unique coefficients
in the Schur function expansions of the plethysms $e_m[e_k]$.
\end{definition}

Computing plethysm coefficients $a_{\lambda,k}$ is a well-known open problem;
see \cite{CarbonaraRemmelYang, Carre, ChenGarsiaRemmel, Howe}.  One of the well-known special cases is when
$k=2$, and is given by the following identity of Littlewood;
see e.g. Macdonald \cite[Exer. I.5.9(a) and I.8.6(c)]{Macdonald}, and also 
Burge \cite{Burge} and Gasharov \cite{Gasharov} for connections with graphical
degree sequences.

\begin{theorem}
\label{Littlewood}
$$
\Psi_2(\xx)=\prod_{i<j} (1 + x_i x_j) = \sum_{m \geq 0} e_m[e_2] =
 \sum_{\text{ shifted }2\text{-families }K} s_{d(K)}. 
$$
\end{theorem}
In other words, for $k=2$, one has 
$$
a_{\lambda,2}=
\begin{cases}
1 & \text{ if } \lambda=d(K) \text{ for some shifted }2\text{-family }K\\
0 & \text{ otherwise.}
\end{cases}
$$

Recall that the {\it Kostka coefficient} $K_{\lambda,\mu}$ is the number of
column-strict tableaux of shape $\lambda$ and content $\mu$, and
gives the expansion 
\begin{equation}
\label{Kostka-expansion}
s_\lambda = \sum_\mu K_{\lambda,\mu} m_\mu.
\end{equation}
It is easily seen and well-known that the $K_{\lambda,\mu}$ are
unitriangular: $K_{\lambda,\lambda}=1$ and
$K_{\lambda,\mu} \neq 0$ if and only if $\lambda \text{ majorizes }\mu$.
Thus knowledge of the plethysm coefficients $a_{\lambda,k}$ determines
the numbers $c_{\mu,k}$ via
\begin{equation}
c_{\mu,k} = \sum_\lambda a_{\lambda,k} K_{\lambda,\mu}.
\end{equation}
In particular, 
$c_{\mu,k} \neq 0$ if and only if there exists some $\lambda$ which has
$a_{\mu,k} \neq 0$ and majorizes $\mu$.

\begin{definition} \rm \ 
Motivated by the form of Theorem~\ref{Littlewood}, we define for each $k \geq 0$
another symmetric function $\Phi_k(\xx)$, and give names to its coefficients
$a'_{\lambda,k}, c'_{\mu,k}$ when expanded in the monomial and Schur function bases:
$$
\begin{aligned}
\Phi_k(\xx) &:= \sum_{\text{ shifted }k\text{-families }K} s_{d(K)} \\
            &= \sum_{\lambda: |\lambda| \equiv 0 \mod k} a'_{\lambda,k}  s_\lambda \\
            &= \sum_{\mu: |\mu| \equiv 0 \mod k} c'_{\mu,k} m_\mu.
\end{aligned}
$$
Note that, by the above definition, 
$a'_{\lambda,k}$ is the number of shifted $k$-families with $d(K)=\lambda$.
\end{definition}

\begin{proposition}
\label{small-k-coincidences}
The symmetric functions $\Psi_k, \Phi_k$ have the same monomial
support, that is, 
$$
c_{\mu,k} = 0 \text{ if and only if } c'_{\mu,k}=0.
$$
Also, one has for $k=0,1,2$ the following three (equivalent) equalities
\begin{equation}
\begin{aligned}
\Psi_k(\xx)   &=\Phi_k(\xx)\\
a_{\lambda,k} &=a'_{\lambda,k} \text{ for all $\lambda$}\\
c_{\mu,k}     &=c'_{\mu,k} \text{ for all $\mu$ }.
\end{aligned}
\end{equation}
\end{proposition}

\begin{proof}
The assertion about monomial supports is a rephrasing of 
Proposition~\ref{unsatisfactory-proposition}, using the above-stated facts about 
Kostka numbers.  The assertion about the equalities
is somewhat trivial for $k=0,1$, and is Littlewood's identity (Theorem~\ref{Littlewood}) for $k=2$.
\end{proof}

The goal of this section is to explore further the link between
$\Psi_k(\xx)$ and $\Phi_k(\xx)$, that is, between the plethysm problem and 
degree sequences of shifted families, when $k \geq 3$.

\subsection{Symmetries}
\label{symmetries-subsection}
There are two obvious symmetries of $k$-families on vertex set $[n]$.  These lead to
 symmetries of the Schur expansion coefficients $c_{\lambda,k}, c'_{\lambda,k}$ for the
symmetric functions $\Psi_k(\xx), \Phi_k(\xx)$ when one works in a finite variable set $x_1,\ldots,x_n$,
that is, by setting $x_{n+1}=x_{n+2}= \cdots =0$.  Note that when working in this finite variable set,
one has these interpretations:
$$
\Psi_k(\xx) := \sum_{k\text{-families }K\text{ on }[n]} \xx^{d(K)} 
 := \prod_{S \in \binom{[n]}{k}} (1 + \xx_S) \\
 := \sum_{m \geq 0} e_m[e_k(x_1,\ldots,x_n)] \\
$$
and
$$
\Phi_k(\xx) := \sum_{\substack{\text{ shifted }k\text{-families }\\K\text{ on }[n]}} s_{d(K)}(x_1,\ldots,x_n).
$$

We will use freely two basic facts about symmetric functions and Schur functions in finite
variable sets.  Recall that $\ell(\lambda)$ denotes the {\it length} or number of parts
in a partition $\lambda$.

\begin{proposition}
The symmetric functions in $n$ variables have as a $(\ZZ-)$basis the
Schur functions 
$$
\{ s_\lambda(x_1,\ldots,x_n) \}_{\ell(\lambda) \leq n},
$$ 
while $s_\lambda(x_1,\ldots,x_n)=0$ if $\ell(\lambda) > n$.

In particular, one has the following consequence:  if
$f(x_1,x_2,\ldots)$ is a symmetric function in the 
infinite variable set $x_1,x_2,\ldots$ with
(unique) expansion
$$
f=\sum_\lambda  a_\lambda s_\lambda,
$$
then the coefficients 
$a_\lambda$ for $\ell(\lambda) \leq n$ are determined by the
unique expansion of the specialization to $x_1,\ldots,x_n$:
$$
f_\lambda(x_1,\ldots,x_n,0,0,\ldots) 
  = \sum_{\lambda: \ell(\lambda) \leq n} a_\lambda s_\lambda(x_1,\ldots,x_n).
$$
\end{proposition}

\begin{proposition}\cite[Exercise 7.41]{Stanley-EC2}
\label{Schur-at-reciprocals}
Assume $\ell(\lambda) \leq n$ and $\lambda_1 \leq  N$,
so that the Ferrers diagram for $\lambda$ fits inside an $n \times N$ rectangle 
$R=\underbrace{(N,\ldots,N)}_{n\text{ times}}$.
Then 
$$
(x_1 \cdots x_n)^N s_\lambda(x_1^{-1},\ldots,x_n^{-1}) = s_{R \setminus \lambda} (x_1,\ldots,x_n)
$$
where $R \setminus \lambda$ denotes the Ferrers diagram obtained by removing $\lambda$ from the northwest corner of
$R$ and then rotating 180 degrees.
\end{proposition}

We come now to the first symmetry.  There is an involution
on the collection of all $k$-families on $[n]$, which maps $K \mapsto \binom{[n]}{k} \setminus K$.

\begin{proposition}
\label{symmetry-one}
Fix $n$ and $k$, and let $R$ be an $n \times \binom{n-1}{k-1}$ rectangle.  
Then whenever $\ell(\lambda) \leq n$, one has
\begin{enumerate}
\item[(i)] $a_{\lambda,k} = 0$ unless $\lambda_1 \leq \binom{n-1}{k-1}$, in which case $a_{R \setminus \lambda, k}=a_{\lambda,k}$.
\item[(ii)] $a'_{\lambda,k} = 0$ unless $\lambda_1 \leq \binom{n-1}{k-1}$, in which case $a'_{R \setminus \lambda, k}=a'_{\lambda,k}$.
\end{enumerate}

In particular, the plethysm coefficients 
$$
\left\{ a_{\lambda,k}: \ell(\lambda) \leq n \text{ and }
|\lambda| \leq \frac{n}{2}\binom{n-1}{k-1} 
% \left( =\frac{k}{2}\binom{n}{k} \right) 
\right\}
$$
already determine the rest of the $\{ a_{\lambda,k}: \ell(\lambda) \leq n\}$
(and similarly for $a'_{\lambda,k}$).
\end{proposition}

\begin{proof}
The fact that $a_{\lambda,k}=a'_{\lambda,k}=0$ unless $\lambda_1 \leq \binom{n-1}{k-1}$ follows
because the degree sequence $K$ of any $k$-family on vertex set $[n]$ (whether shifted or not)
is bounded above by 
$$
R:=\underbrace{\left(\binom{n-1}{k-1},\ldots,\binom{n-1}{k-1}\right)}_{n\text{ times}}.
$$
Hence if one had a Schur function $s_\lambda$ with $\lambda_1 > \binom{n-1}{k-1}$
in the Schur expansion of $\Psi_k(x_1,\ldots,x_n)$, the leading term $\xx^\lambda$ which occurs
in the monomial expansion of $s_\lambda$ would lead to a contradiction.

Note that the involution $K \mapsto \binom{[n]}{k} \setminus K$ has the
property that 
$$
d\left(\binom{[n]}{k} \setminus K \right)=R - d(K) 
= \left(\binom{n-1}{k-1},\ldots,\binom{n-1}{k-1}\right)-d(K)
$$
where one should be careful to note that $R - d(K)$ is an ordered
degree sequence, in weakly {\it increasing} order, not decreasing order,
so that it is the {\it reverse} of the partition $R \setminus d(K)$.  This implies
$$
\begin{aligned}
(x_1 \cdots x_n)^{\binom{n-1}{k-1}}\Psi_k(x_1^{-1},\ldots,x_n^{-1})
  = \Psi_k(x_1,\ldots,x_n)
\end{aligned}
$$
and using Proposition~\ref{Schur-at-reciprocals} then gives the remaining assertion in (i).

For (ii), note that if one follows this symmetry $K \mapsto \binom{[n]}{k} \setminus K$ by
the map which reverses the vertex labels $i \mapsto n+1-i$ in $[n]$, the composite
is an involution $K \mapsto K^c$ on the collection of
all {\it shifted} $k$-families.  This composite involution satisfies
$d(K^c) = R \setminus d(K)$,
which shows the remaining assertion in (ii).
\end{proof}

There is a second involution that sends $k$-families on $[n]$ to $(n-k)$-families on $[n]$,
by mapping $K\mapsto \{[n] \setminus S: S \in K\}$. 

\begin{proposition}
\label{symmetry-two}
Fix $k$ and $n$.  Assume that
$\ell(\lambda) \leq n$ and $|\lambda| \equiv 0 \mod k$, with $m:=\frac{|\lambda|}{k}$.  
Let $M$ be an $n \times m$ rectangle.  Then
\begin{enumerate}
\item[(i)] $a_{\lambda,k} = 0$ unless $\lambda_1 \leq m$, in which case $a_{M \setminus \lambda, k}=a_{\lambda, n-k}$.
\item[(ii)] $a'_{\lambda,k} = 0$ unless $\lambda_1 \leq m$, in which case $a'_{M \setminus \lambda, k}=a'_{\lambda, n-k}$.
\end{enumerate}

Thus the plethysm coefficients $\{ a_{\lambda,k}: \ell(\lambda) \leq n \}$ for a fixed $k$ determine
the plethysm coefficients $\{a_{\lambda,n-k}:\ell(\lambda) \leq n\}$.
In particular, Proposition~\ref{small-k-coincidences} determines all the 
$a_{\lambda,k}$ both for $k \leq 2$, and for $\ell(\lambda) \leq k+2$.
\end{proposition}

\begin{proof}
The fact that $a_{\lambda,k}=a'_{\lambda,k}=0$ unless $\lambda_1 \leq m$ follows
because a $k$-family $K$ on vertex set $[n]$ of cardinality $|K|=m$ (meaning there are $m$ sets in $K$) will have
degree sequence $d(K)$ bounded above by
$$
M:=\underbrace{(m,\ldots,m)}_{n\text{ times}}.
$$
Hence if one had a Schur function $s_\lambda$ with $\lambda_1 > m$
in the Schur expansion of $\Psi_k(x_1,\ldots,x_n)$, the term leading $\xx^\lambda$ which occurs
in the monomial expansion of $s_\lambda$ would lead to a contradiction.

For (i), it is not hard to check (e.g. using the involution 
$K\mapsto \{[n] \setminus S: S \in K\}$) that
$$
(x_1 \cdots x_n)^{\binom{n}{k}} \Psi_k(x_1^{-1},\ldots,x_n^{-1}) = \Psi_{n-k}(x_1,\ldots,x_n),
$$
and the assertion then follows from Proposition~\ref{Schur-at-reciprocals}.

For (ii), note that if one follows this symmetry $K\mapsto \{[n] \setminus S: S \in K\}$
by the map  which reverses the vertex labels $i \mapsto n+1-i$ in $[n]$,
one obtains an involution $K \mapsto K^*$ on the collection of all {\it shifted} families. 
This composite involution satisfies $d(K^*) = M \setminus d(K)$,
which shows the remaining assertion in (ii)
\end{proof}

\subsection{Highest weight vectors}

The goal of this section is the following result relating
the Schur expansion coefficients $a_{\lambda,k}, a'_{\lambda,k}$ for $\Psi_k(\xx), \Phi_k(\xx)$.

\begin{theorem}
\label{shifteds-as-lower-bound}
For all $k$ and $\lambda$, one has 
$$
a_{\lambda, k} \geq a'_{\lambda, k}.
$$
Furthermore, equality holds when either $k=2$ or $\ell(\lambda) \leq k+2$.
\end{theorem}

Since the statement about equality follows from Proposition~\ref{small-k-coincidences} and 
Proposition~\ref{symmetry-two}, we must only show the stated inequality.

For this, we review a standard $GL_n$-representation interpretation for the plethysm $e_m[e_k]$;
see \cite[Appendix A.7, Example]{Macdonald} and \cite[Chapter 7, Appendix 2]{Stanley-EC2}.
Let $V=\CC^n$ with a chosen $\CC$-basis $e_1,\ldots,e_n$.  Then $G=GL_n(\CC)$ acts on $V$,
and if one chooses as a maximal torus $T$ the diagonal matrices in $G$,
the typical element $x:=\diag(x_1,\ldots,x_n)$ in $T$ has $e_i$ as an eigenvector with
eigenvalue $x_i$.  In other words, the $\{ e_i \}$ form a basis of {\it weight} vectors for $T$.
The {\it character} of any (polynomial) $G$-representation $U$ is defined to be the symmetric function 
$\character(U)$ in the variables $x_1,\ldots,x_n$ which is the trace of $x$ acting on $U$, 
or the sum of the eights/eigenvalues of $x$ in any basis of $T$-weight vectors for $U$.  Thus for 
$V=\CC^n$ one has $\character(V)=x_1+\ldots+x_n$.

The $k^{th}$ exterior power $\bigwedge^k V$ inherits a $\CC$-basis of monomial decomposable
wedges, indexed by $k$-sets $S=\{i_1 < \cdots < i_k\}$, and defined by
$$
e_S:=e_{i_1} \wedge \cdots \wedge e_{i_k}.
$$
Furthermore, each $e_S$ is a $T$-weight vector with weight $\xx_S:=x_{i_1} \cdots x_{i_k}$, and
hence 
$$
\character\left( \bigwedge^k V \right)=\sum_{S \in \binom{[n]}{k}} \xx_S = e_k(x_1,\ldots,x_n).
$$

The $m^{th}$ exterior power $\bigwedge^m \left(\bigwedge^k V \right)$ similarly inherits
a $\CC$-basis of monomial decomposable wedges, indexed by $k$-families $K=\{S_1,\ldots,S_m\}$
of cardinality $|K|=m$, and defined by
$$
E_K := e_{S_1} \bigwedge \cdots \bigwedge e_{S_m}.
$$
We will assume that the $k$-sets $S_1,\ldots,S_m$ in $K$ are always listed in some
fixed linear ordering (such as lexicographic order), for the sake of definiteness in
writing down $E_K$;  the chosen order will only affect $E_K$ up to scaling by $\pm 1$.
Furthermore, each $E_K$ is a $T$-weight vector with weight $\prod_{S \in K} \xx_S = \xx^{d(K)}$,
and hence
$$
\character\left(\bigwedge^m \left(\bigwedge^k V \right) \right)
  = \sum_{\substack{ k\text{-families }K\text{ on }[n]\\ \text{ of cardinality }m}} \xx^{d(K)} = e_m[e_k].
$$

As this $GL_n(\CC)$-representation $\bigwedge^m \left(\bigwedge^k V \right)$ is polynomial, 
it can be written as a direct sum of the irreducible 
polynomial representations $W^\lambda$, which
are parametrized by partitions $\lambda$ with $\ell(\lambda) \leq n$.  The irreducible
$W^\lambda$ has $\character(W^\lambda) = s_\lambda(x_1,\ldots,x_n)$.
The plethysm coefficient $a_{\lambda,k}$ is then 
the multiplicity of the irreducible $W^\lambda$ in the decomposition of 
$\bigwedge^m \left(\bigwedge^k V \right)$.  Because $\bigwedge^m \left(\bigwedge^k V \right)$ is
homogeneous of degree $km$, the $W^\lambda$ which appear in the decomposition will have
$|\lambda|=km$.

To bring in the shifted families, it is convenient to use highest weight theory for the associated
lie algebra $\lieg=\lieg\ell_n(\CC)$ of all $n \times n$ matrices over $\CC$.  An $n \times n$ matrix $A$
in $\lieg$ acts on $V=\CC^n$ in the usual way.  Once one knows the action of an element $A \in \lieg$ on
a space $U$, then it acts on $\bigwedge^k U$ by a Leibniz rule:  
$$
A\left( u_1 \wedge \cdots \wedge u_k \right) =
 \sum_{i=1}^k u_1 \wedge \cdots \wedge u_{i-1} \wedge  (Au_i) \wedge u_{i+1} \wedge \cdots \wedge u_k.
$$

The decomposition for a polynomial $G$-representation $U$
into irreducibles is determined by the subspace of highest weight vectors in $U$.  Specifically, if
one chooses as a nilpotent subalgebra $\nn_+$ the set of all strictly upper triangular matrices in $\lieg$,
then the subspace $U^+$ of $U$ annihilated by all of $\nn_+$ is called the space of {\it highest weight vectors}.
It turns out that $U^+$ will always have a basis $\{u_1,\ldots,u_p\}$ 
of $T$-weight vectors $u_i$, each having a weight which is
{\it dominant}, that is, of the form $\xx^{\lambda^{(i)}}$ for some partition $\lambda^{(i)}$.  The theory
tells us that then the irreducible decomposition of $U$ is 
$$
U \cong \bigoplus_{i=1}^p W^{\lambda^{(i)}}.
$$
Consequently, $\character(U) = \sum_{i=1}^p s_{\lambda^{(i)}}.$

Based on the previous discussion, Theorem~\ref{shifteds-as-lower-bound} follows immediately
from the next proposition\footnote{Since one can alternately define the highest weight vectors in a 
$GL_n(\CC)$-representation as those fixed by the elements of the {\it Borel subgroup} of upper triangular
invertible matrices, Proposition~\ref{shifteds-are-highest-weights} can be viewed as an exterior algebra analogue
to the combinatorial description of Borel-fixed monomial ideals in the symmetric algebra of $V$;
see e.g. \cite[Proposition 2.3]{MillerSturmfels}}.

\begin{proposition}
\label{shifteds-are-highest-weights}
Among the basis of monomial weight vectors $E_K$ for $\bigwedge^m \left(\bigwedge^k V \right)$ indexed
by the $k$-families  $K$ on $[n]$, those which are highest weight vectors are
exactly the ones indexed by shifted $k$-families $K$.
\end{proposition}
\begin{proof}
The subalgebra $\nn_+$ has a $\CC$-basis of elementary matrices $\{ A_{i,j}: 1 \leq i < j \leq n\}$,
where $A_{i,j}$ has one non-zero entry, equal to $1$, in column $j$, row $i$.  In other words,
$A_{i,j}$ acts on $V$ by killing all basis vectors $e_r$ except for $e_j$, which it sends to $e_i$.
By the Leibniz rule, one checks that $A_{i,j}$ acts on the monomial basis vectors 
$e_S$ for
$\bigwedge^k V$ as follows:  $A_{i,j}$ kills $e_S$ unless $j \in S$ and $i \not\in S$, in which
case $A_{i,j}(e_S)=\pm e_{S'}$ where $S':= S \setminus \{j\} \cup \{i\}$.  Note that $S'$ is
a set strictly lower in the componentwise ordering than $S$.

Given a $k$-family $K=\{S_1,\ldots,S_m\}$, by the Leibniz rule one has
$$
A_{i,j}(E_K) =
 \sum_{\ell=1}^m e_{S_1} \bigwedge \cdots \bigwedge A_{i,j}(e_{S_\ell}) \bigwedge \cdots \bigwedge e_{S_m}.
$$
If $K$ is a shifted $k$-family, we claim this sum vanishes term-by-term: either $A_{i,j}(e_{S_\ell})=0$,
or $A_{i,j}(e_{S_\ell})=\pm e_{S'}$ for some $S'$ lower in the componentwise order than $S$, 
so that $e_{S'}$ coincides with another element in the wedge and the term still vanishes.

Conversely, if $K$ is {\it not} a shifted $k$-family, there must exists at least one set $S$ in $K$
and some pair of indices $i < j$ for which 
\begin{enumerate}
\item[$\bullet$]
$j \in S$,
\item[$\bullet$]
$i \not\in S$, and
\item[$\bullet$]
$S':= S\setminus \{j\} \cup \{i\}$ is {\it not} in $K$.
\end{enumerate}
In this case, assume by re-indexing that $S_1,\ldots,S_{r}$ are the sets which have this
property (so $S'_1,\ldots,S'_r$ are well-defined), 
and $S_{r+1},\ldots,S_{m-1},S_m$ are the ones that do not.  Then the above calculation shows 
$$
A_{i,j}(E_K) = 
  \sum_{\ell=1}^r  \pm e_{S_1} \bigwedge \cdots \bigwedge e_{S'_\ell} \bigwedge \cdots \bigwedge e_{S_m}
$$
and one can check that each term in this sum is (up to $\pm 1$) a {\it different} one
of the monomial basis vectors $E_{K'}$:  if two such terms indexed by $\ell, \ell'$
were to coincide, their corresponding sets $S'_\ell, S'_{\ell'}$ would need to coincide,
forcing the sets $S_\ell, S_{\ell'}$ to coincide,  i.e. $\ell=\ell'$.
So $A_{i,j}(E_K)$ does not vanish.
\end{proof}

\begin{openproblem}
What more can one say about the (Schur-positive) difference $\Psi_k(\xx) - \Phi_k(\xx)$?
\end{openproblem}

We offer a conjecture in this direction, which considers the various homogeneous
components of this difference.  If true, it suggests that when computing a
plethysm $e_m[e_k]$, not only are the shifted $k$-families of size $m$ relevant for the
{\it top} of the expansion, but those of size $i < m$ seem relevant for the rest of the
expansion.

\begin{definition} \rm \ 
Let  $\Phi_{k,m}(\xx)$ be the homogeneous component of $\Phi_k$ having
degree $km$ in the variables $x_i$, that is
$$
\Phi_{k,m}(\xx) := \sum_{\substack{\text{ shifted }k\text{-families }\\K\text{ s.t. }|K| = m}} s_{d(K)}.
$$
\noindent
The analogous homogeneous component of $\Psi_k(\xx)$ is
the plethysm $e_m[e_k]$.  

Define a system $\Upsilon_{k,m}(\xx)$ as follows:

\begin{eqnarray*}
\Upsilon_{k,1}(\xx) &:=  & e_1[e_k] - \Phi_{k,1}(\xx) =0.\\
\Upsilon_{k,2}(\xx) & := & e_2[e_k] - \Phi_{k,2}(\xx) \\
\Upsilon_{k,3}(\xx) & := & e_3[e_k] - \Phi_{k,3}(\xx) - \Upsilon_{k,2}(\xx)  \Phi_{k,1}(\xx) \\
\Upsilon_{k,4}(\xx) & := & e_4[e_k] - \Phi_{k,4}(\xx) - \Upsilon_{k,3}(\xx)  \Phi_{k,1}(\xx)
                                                      -\Upsilon_{k,2}(\xx)  \Phi_{k,2}(\xx) \\
\vdots & & \\
\Upsilon_{k,m}(\xx) &:= & e_m[e_k] - \Phi_{k,m}(\xx) - \sum_{i = 1}^{m-2} (\Upsilon_{k,m-i}(\xx)  \Phi_{k,i}(\xx)).
\end{eqnarray*}
\end{definition}

\begin{conjecture}

$(-1)^m \Upsilon_{k,m}(\xx)$ is Schur-positive.

\end{conjecture}

For $m=1$, this conjecture is trivial. 

For $m=2$, it is nearly trivial, and is
implied by Theorem~\ref{shifteds-as-lower-bound}.

For $m=3$, it has been checked using explicit expansions
of $e_3[e_k]$, such as the one given by Chen, Garsia and Remmel \cite{ChenGarsiaRemmel}.  

For $m=4$, although in principle one might be able to check it using
the explicit expansions of $e_4[e_k]$ given by Foulkes \cite{Foulkes}
and Howe \cite{Howe}, in practice the calculations are unpleasant
enough that we have not done them.  

We have also checked (using Stembridge's Maple package for symmetric functions
{\tt SF}) that the conjecture holds for these values:
$$
\begin{aligned}
m = 4  ,  & \text{ and } k \leq 7\\
m = 5  ,  & \text{ and } k \leq 6\\
m = 6  ,  & \text{ and } k \leq 4\\
m = 7  ,  & \text{ and } k \leq 4\\
m = 8  ,  & \text{ and } k \leq 3.
\end{aligned}
$$

\section*{Acknowledgments}
The authors thank Andrew Crites, Pedro Felzenszwalb and two anonymous referees for
helpful edits, comments and suggestions.

\end{document}